\documentclass[10pt,a4paper]{amsart}
\usepackage{graphicx}
\usepackage{xcolor}
\usepackage{latexsym}
\usepackage{verbatim}
\usepackage{amsmath}
\usepackage{amsthm}
\usepackage{amssymb}
\usepackage{hyperref}
\numberwithin{equation}{section}
\newtheorem{theorem}{Theorem}[section]

\newtheorem{remark}[theorem]{Remark}

\title[A Numerical Scheme for the Quantum Boltzmann Equation Efficient in the Fluid Regime]{A Numerical Scheme for the Quantum Boltzmann Equation Efficient in the Fluid Regime}\thanks{This work was partially supported by NSF grant DMS-0608720 and NSF FRG grant DMS-0757285. FF was supported by the ERC Starting Grant Project NuSiKiMo. SJ was also supported by a Van Vleck Distinguished Research Prize and a Vilas Associate Award from the University of Wisconsin-Madison.}

\def\signff{\bigskip\bigskip\hspace{80mm}
\vbox{{\sc Francis Filbet\par\vspace{3mm}
Universit\'e de Lyon,\par
Universit\'e Lyon I, CNRS \par
UMR 5208, Institut Camille Jordan \par
43, Boulevard du 11 Novembre 1918\par
69622 Villeurbanne cedex, FRANCE\par\vspace{3mm}
e-mail:} filbet@math.univ-lyon1.fr }}

\def\signjh{\bigskip\bigskip\hspace{80mm}
\vbox{{\sc Jingwei Hu\par\vspace{3mm}
Department of Mathematics, \par 
University of Wisconsin-Madison, \par
480 Lincoln Drive, Madison,  \par
WI 53706, USA \par\vspace{3mm}
email:} hu@math.wisc.edu}}

\def\signsj{\bigskip\bigskip\hspace{80mm}
\vbox{{\sc Shi Jin\par\vspace{3mm}
Department of Mathematics, \par 
University of Wisconsin-Madison,
\par
480 Lincoln Drive, Madison,  \par
WI 53706, USA \par\vspace{3mm}
email:} jin@math.wisc.edu}}

\author{Francis Filbet , Jingwei Hu \and Shi Jin}

\begin{document}

\maketitle

\begin{abstract}
Numerically solving the Boltzmann kinetic equations with the small Knudsen number is challenging due to the stiff nonlinear collision term. A class of asymptotic preserving schemes was introduced in \cite{FJ10} to handle this kind of problems. The idea is to penalize the stiff collision term by a BGK type operator. This method, however, encounters its own difficulty when applied to the quantum Boltzmann equation. To define the quantum Maxwellian (Bose-Einstein or Fermi-Dirac distribution) at each time step and every mesh point, one has to invert a nonlinear equation that connects the macroscopic quantity fugacity with density and internal energy. Setting a good initial guess for the iterative method is troublesome in most cases because of the complexity of the quantum functions (Bose-Einstein or Fermi-Dirac function). In this paper, we propose to penalize the quantum collision term by
a `classical' BGK operator instead of the quantum one. This is based on the observation that the classical Maxwellian, with the temperature replaced by the internal energy, has the same first five moments as the quantum Maxwellian. The scheme so designed avoids the aforementioned difficulty, and one can show that the density distribution is still driven toward the quantum equilibrium. Numerical results are present to illustrate the efficiency of the new scheme in both the hydrodynamic and kinetic regimes. We also develop a spectral method for the quantum collision operator.
\end{abstract}

\section{Introduction}
The quantum Boltzmann equation (QBE), also known as the Uehling-Uhlenbeck equation, describes the behaviors of a dilute quantum gas. It was first formulated by Nordheim \cite{Nordheim28} and Uehling and Uhlenbeck \cite{UU33} from the classical Boltzmann equation by heuristic arguments. Here we mainly consider two kinds of quantum gases: the Bose gas and the Fermi gas. The Bose gas is composed of Bosons, which have an integer value of spin, and obey the Bose-Einstein statistics. The Fermi gas is composed of Fermions, which have half-integer spins and obey the Fermi-Dirac statistics. 

Let $f(t,x,v)\geq 0$ be the phase space distribution function depending on time $t$, position $x$ and particle velocity $v$, then the quantum Boltzmann equation reads:
\begin{equation} \label{QBE}
\frac{\partial f}{\partial t}+v\cdot \nabla _x f=\frac{1}{\epsilon}\mathcal{Q}_q(f),\   \    \  x\in \Omega \subset \mathbb{R}^{d_x},\ v\in \mathbb{R}^{d_v}.
\end{equation}
Here $\epsilon$ is the Knudsen number which measures the degree of rarefaction of a gas. It is the ratio between the mean free path and the typical length scale. The quantum collision operator $\mathcal{Q}_q$ is
\begin{equation} \label{Qq}
 \mathcal{Q}_q(f)(v)=\int _{\mathbb{R}^{d_v}} \int_{\mathbb{S}^{d_v-1}}B(v-v_*,\omega)\left [f'f_*'(1\pm \theta_0 f)(1\pm \theta_0 f_*)-ff_*(1\pm \theta_0 f')(1\pm \theta_0f_*') \right ]d\omega dv_*
\end{equation}
where $\theta_0=\hbar^{d_v}$, $\hbar$ is the rescaled Planck constant. In this paper, the upper sign will always correspond to the Bose gas while the lower sign to the Fermi gas. For the Fermi gas,  we also need $f\leq \frac{1}{\theta_0}$ by the Pauli exclusion principle. $f$, $f_*$, $f'$ and $f_*'$ are the shorthand notations for $f(t,x,v)$, $f(t,x,v_*)$, $f(t,x,v')$ and $f(t,x,v_*')$ respectively. $(v,v_*)$ and $(v',v_*')$ are the velocities before and after collision. They are related by the following parametrization:
\begin{eqnarray}
\left\{
\begin{array}{l}
\displaystyle v'=\frac{v+v_*}{2}+\frac{|v-v_*|}{2}\omega, 
\\
\\
\displaystyle v_*'=\frac{v+v_*}{2}-\frac{|v-v_*|}{2}\omega,
\end{array}\right.
\end{eqnarray}
where $\omega$ is the unit vector along $v'-v_*'$. The collision kernel $B$ is a nonnegative function that only depends on $|v-v_*|$ and $\cos \theta$ ($\theta$ is the angle between $\omega$ and $v-v_*$). In the Variable Hard Sphere (VHS) model, it is given by
\begin{equation}
B(v-v_*,\omega)=C_\gamma |v-v_*|^{\gamma}
\end{equation}
where $C_\gamma$ is a positive constant. $\gamma=0$ corresponds to the Maxwellian molecules, $\gamma=1$ is the hard sphere model.

When the Knudsen number $\epsilon$ is small, the right hand side of equation (\ref{QBE}) becomes stiff and explicit schemes are subject to severe stability constraints. Implicit schemes allow larger time step, but new difficulty arises in seeking the numerical solution of a fully nonlinear problem at each time step. Ideally, one wants an implicit scheme allowing large time steps and can be inverted easily. In \cite{FJ10}, for the classical Boltzmann equation, Filbet and Jin proposed to penalize the nonlinear collision operator $\mathcal{Q}_c$ by a BGK operator:
\begin{equation} \label{pen}
\mathcal{Q}_c=[\mathcal{Q}_c-\lambda(\mathcal{M}_c-f)]+\lambda[\mathcal{M}_c-f]
\end{equation}
where $\lambda$ is a constant that depends on the spectral radius of the linearized collision operator of $\mathcal{Q}_c$ around the local (classical) Maxwellian $\mathcal{M}_c$. Now the term in the first bracket of the right hand side of (\ref{pen}) is less stiff than the second one and can be treated explicitly. The term in the second bracket will be discretized implicitly. Using the conservation property of the BGK operator, this implicit term can actually be solved explicitly. Thus they arrive at a scheme which is uniformly stable in $\epsilon$, with an implicit source term that can be inverted explicitly. Furthermore, under certain conditions, one could show that this type of schemes has the following property: the distance between $f$ and the Maxwellian will be $O(\epsilon)$ after several time steps, no matter what the initial condition is. This guarantees the capturing of the fluid dynamic limit even if the time step is larger than the mean free time.

Back to the quantum Boltzmann equation (\ref{QBE}), a natural way to generalize the above idea is to penalize $\mathcal{Q}_q$ with the quantum BGK operator $\mathcal{M}_q-f$. This means we have to invert a nonlinear algebraic system that
contains the unknown  quantum Maxwellian $\mathcal{M}_q$ (Bose-Einstein or Fermi-Dirac distribution) for every time step. As mentioned in \cite{hu_KFVS}, this is not a trivial task compared to the classical case. Specifically,
one has to invert a nonlinear 2 by 2 system (can be reduced to one nonlinear equation) to obtain the macroscopic quantities, temperature and fugacity. Due to the complexity of the quantum distribution  functions (Bose-Einstein or Fermi-Dirac function), it is really a delicate issue to set a good initial guess for an iterative method such as the Newton method to converge.

In this work we propose a new scheme for the quantum Boltzmann equation. Our idea is based on the observation that the classical Maxwellian, with the temperature replaced by the (quantum) internal energy, has the same first five moments as the quantum Maxwellian. This observation was used in \cite{hu_KFVS} to derive a `classical' kinetic scheme for the quantum hydrodynamical equations. Therefore, we just penalize the quantum collision operator $\mathcal{Q}_q$ by a `classical' BGK operator, thus avoid the aforementioned difficulty. At the same time, we have to sacrifice a little bit on the asymptotic property. Later we will prove that for the quantum BGK equation, the so obtained $f$ satisfies:
\begin{equation} \label{property}
f^n-\mathcal{M}_q^n=O(\Delta t) \  \  \text{ for some }n> N,  \  \text{any initial data } f^0,
\end{equation} 
i.e. $f$ will converge to the quantum Maxwellian beyond the initial layer with an error of $O(\Delta t)$.

Another numerical issue is how to evaluate the quantum collision operator $\mathcal{Q}_q$. In fact (\ref{Qq}) can be simplified as
\begin{equation} \label{Qqq}
 \mathcal{Q}_q(f)(v)=\int _{\mathbb{R}^{d_v}} \int_{\mathbb{S}^{d_v-1}}B(v-v_*,\omega)\left [f'f_*'(1\pm \theta_0 f\pm \theta_0 f_*)-ff_*(1\pm \theta_0 f'\pm \theta_0f_*') \right ]d\omega dv_*
\end{equation}
so $\mathcal{Q}_q$ is indeed a cubic operator. Almost all the existing fast algorithms are designed for the classical Boltzmann operator based on its quadratic structure. Here we will give a spectral method for the approximation of $\mathcal{Q}_q$. As far as we know, this is the first time to compute the full quantum Boltzmann collision operator with the spectral accuracy.

The rest of the paper is organized as follows. In the next section, we give a brief introduction to the quantum Boltzmann equation: the basic properties, the quantum Maxwellians and the hydrodynamic limits. In section 3,  we present the details of computing the quantum collision operator by the spectral method as well as the numerical accuracy. Our new scheme to capture  the hydrodynamic regime is given  in section 4.  In section 5, the proposed schemes are tested on the 1-D shock tube problem of the quantum gas for different Knudsen number $\epsilon$ ranging from fluid regime to kinetic regime. The behaviors of the Bose gas and the Fermi gas in both the classical regime and quantum regime are included. Finally some concluding remarks are given in section 6.

\section{The Quantum Boltzmann Equation and its Hydrodynamic Limits}
In this section we review some basic facts about the quantum Boltzmann equation (\ref{QBE}). 
\begin{itemize}
\item At the formal level, $\mathcal{Q}_q$ conserves mass, momentum and energy.
\begin{equation}
\int_{\mathbb{R}^{d_v}}\mathcal{Q}_q(f)dv=\int_{\mathbb{R}^{d_v}}\mathcal{Q}_q(f)vdv=\int_{\mathbb{R}^{d_v}}\mathcal{Q}_q(f)|v|^2dv=0.
\end{equation}
\item If $f$ is a solution of QBE (\ref{QBE}), the following local conservation laws hold:
\begin{eqnarray}
\left\{
\begin{array}{l}
\displaystyle\frac{\partial}{\partial t} \int_{\mathbb{R}^{d_v}}fdv+\nabla_x\cdot \int_{\mathbb{R}^{d_v}}vfdv=0, 
\\
\\
\displaystyle\frac{\partial}{\partial t} \int_{\mathbb{R}^{d_v}}vfdv+\nabla_x\cdot \int_{\mathbb{R}^{d_v}}v\otimes vfdv=0,
\\
\\
\displaystyle\frac{\partial}{\partial t} \int_{\mathbb{R}^{d_v}}\frac{1}{2}|v|^2fdv+\nabla_x\cdot \int_{\mathbb{R}^{d_v}}v\frac{1}{2}|v|^2fdv=0.
\end{array}\right.
\label{conserv}
\end{eqnarray}
Define the macroscopic quantities: density $\rho$, macroscopic velocity $u$, specific internal energy $e$ as
\begin{eqnarray}
\label{den}
\rho=\int_{\mathbb{R}^{d_v}}fdv, \quad \rho \, u=\int_{\mathbb{R}^{d_v}}vfdv, \quad \rho e=\int_{\mathbb{R}^{d_v}}\frac{1}{2}|v-u|^2fdv
\end{eqnarray}
and stress tensor  $\mathbb{P}$ and heat flux $q$
\begin{eqnarray} 
 \mathbb{P}=\int_{\mathbb{R}^{d_v}}(v-u)\otimes(v-u)fdv, 
\label{stress}\quad q=\int_{\mathbb{R}^{d_v}}\frac{1}{2}(v-u)|v-u|^2fdv, 
\end{eqnarray}
the above system can then be recast as
\begin{eqnarray}
\left\{
\begin{array}{l}
\displaystyle\frac{\partial\rho}{\partial t} \,+\,\nabla_x \cdot(\rho u)=0, 
\\
\\
\displaystyle\frac{\partial (\rho u)}{\partial t} \,+\,\nabla_x\cdot \left(\rho u\otimes u +\mathbb{P}\right)=0, 
\\
\\
\displaystyle\frac{\partial}{\partial t}\left(\rho e+\frac{1}{2}\rho u^2\right)\,+\,\nabla_x\cdot\left ( \left(\rho e+\frac{1}{2}\rho u^2\right)u+\mathbb{P}u+q\right)=0. 
\end{array}\right.
\label{euler}
\end{eqnarray}

\item $\mathcal{Q}_q$ satisfies Boltzmann's H-Theorem,
\begin{equation}
\int_{\mathbb{R}^{d_v}}\ln \left(\frac{f}{1\pm \theta_0f}\right)\,\mathcal{Q}_q(f)dv\,\leq\, 0,
\end{equation}
moreover,
\begin{equation}
\int _{\mathbb{R}^{d_v}}\ln \left(\frac{f}{1\pm \theta_0 f}\right) \mathcal{Q}_q(f)dv=0 \Longleftrightarrow \mathcal{Q}_q(f)=0  \Longleftrightarrow  f=\mathcal{M}_q,
\end{equation}

where $\mathcal{M}_q$ is the quantum Maxwellian given by
\begin{equation} \label{max}
\mathcal{M}_q=\frac{1}{\theta_0}\frac{1}{z^{-1}e^{\frac{(v-u)^2}{2T}}\mp1},
\end{equation}
where $z$ is the fugacity, $T$ is the temperature (see \cite{hu_KFVS} for more details about the derivation of $\mathcal{M}_q$). This is the well-known Bose-Einstein (`-') and Fermi-Dirac (`+') distributions. 
\end{itemize}

\subsection{Hydrodynamic Limits}
Substituting $\mathcal{M}_q$ into (\ref{den}) (\ref{stress}), the system (\ref{euler}) can be closed, yielding the quantum Euler equations:
\begin{eqnarray}
\label{QEE}
\left\{
\begin{array}{l}
\displaystyle\frac{\partial\rho}{\partial t} \,+\,\nabla_x \cdot(\rho u)=0, 
\\
\\
\displaystyle\frac{\partial(\rho u)}{\partial t} \,+\,\nabla_x\cdot \left(\rho u\otimes u +\frac{2}{d_v}\rho eI \right)=0, 
\\
\\
\displaystyle\frac{\partial}{\partial t}\left(\rho e+\frac{1}{2}\rho u^2\right)\,+\,\nabla_x\cdot \left( \left(\frac{d_v+2}{d_v}\rho e+\frac{1}{2}\rho u^2\right)u\right)=0. 
\end{array}\right.
\end{eqnarray}
With the macroscopic variables $\rho$, $u$ and $e$, they are exactly the same as the classical Euler equations. However, the intrinsic constitutive relation is quite different. $\rho$ and $e$ are connected with $T$ and $z$ (used in the definition of $\mathcal{M}_q$ (\ref{max})) by a nonlinear 2 by 2 system:
\begin{eqnarray}
\label{rhoue}
\left\{
\begin{array}{l}
\displaystyle  \rho=\frac{(2\pi T)^{\frac{d_v}{2}}}{\theta_0}Q_{\frac{d_v}{2}}(z), 
\\
\\
\displaystyle e=\frac{d_v}{2}T\frac{Q_{\frac{d_v+2}{2}}(z)}{Q_{\frac{d_v}{2}}(z)}, 
\end{array}\right.
\end{eqnarray}
where $Q_{\nu}(z)$ denotes the Bose-Einstein 
function $G_{\nu}(z)$ and the Fermi-Dirac function $F_{\nu}(z)$ respectively,
\begin{eqnarray}
G_{\nu}(z)&=&\frac{1}{\Gamma (\nu) }\int _0^\infty  \frac{x^{\nu-1}}{z^{-1}e^x-1}dx,  \  \ 0<z <1, \  \nu>0;  \  z=1,\  \nu>1,  \label{bose}\\
F_{\nu}(z)&=&\frac{1}{\Gamma (\nu) }\int _0^\infty  \frac{x^{\nu-1}}{z^{-1}e^x+1}dx,  \  \ 0<z <\infty, \  \nu>0,  \label{fermi}
\end{eqnarray}
and $\Gamma(\nu)=\int_0^\infty x^{\nu-1}e^{-x}dx$ is the Gamma function. 

The physical range of interest for a Bose gas is $0<z\leq1$, where $z=1$ corresponds to the degenerate case (the onset of Bose-Einstein condensation). For the Fermi gas we don't have such a
restriction and the degenerate case is reached when $z$ is very large. For small $z$ ($0<z<1$), the integrand in (\ref{bose}) and (\ref{fermi}) can be expanded in powers of $z$,
\begin{eqnarray} 
&&G_{\nu}(z)=\displaystyle \sum_{n=1}^{\infty}\frac{z^n}{n^{\nu}}=z+\frac{z^2}{2^{\nu}}+\frac{z^3}{3^{\nu}}+\dots, \label{expan1}\\  
&&F_{\nu}(z)=\displaystyle \sum_{n=1}^{\infty}(-1)^{n+1}\frac{z^n}{n^{\nu}}=z-\frac{z^2}{2^{\nu}}+\frac{z^3}{3^{\nu}}-\dots. \label{expan2}
\end{eqnarray}
Thus, for $z\ll1$, both functions behave like $z$ itself and one recovers
 the classical limit. 
 
On the other hand, the first equation of (\ref{rhoue}) can be written as
\begin{equation}
Q_{\frac{d_v}{2}}(z)=\frac{\rho}{(2\pi T)^{\frac{d_v}{2}}}\theta_0
\end{equation}
where $\frac{\rho}{(2\pi T)^{\frac{d_v}{2}}}$ is just the coefficient of the classical Maxwellian, which should be an $O(1)$ quantity. Now if 
 $\theta_0 \rightarrow 0$, then $Q_{\frac{d_v}{2}}(z) \rightarrow 0$, which means $z\ll 1$ by the monotonicity of the function $Q_{\nu}$. This is consistent with the fact that one gets the classical Boltzmann equation in QBE (\ref{QBE}) by letting $\theta_0 \rightarrow 0$.

The quantum Euler equations (\ref{QEE}) can be derived via the Chapman-Enskog expansion \cite{Cercignani} as the leading order approximation of the quantum Boltzmann equation (\ref{QBE}). By going to the next order, one can also obtain the quantum Navier-Stokes system which differs from their classical counterparts. In particular, the viscosity coefficient and the heat conductivity depend upon both $\rho$ and $e$ \cite{AL97}.

\section{Computing the Quantum Collision Operator $\mathcal{Q}_q$}
In this section, we discuss the approximation of the quantum collision operator $\mathcal{Q}_q$. The method we use is an extension of the spectral method introduced in \cite{MP06,FMP06} for the classical collision operator.

We first write (\ref{Qq}) as
\begin{equation}
\mathcal{Q}_q=\mathcal{Q}_c\pm\theta_0(\mathcal{Q}_1+\mathcal{Q}_2-\mathcal{Q}_3-\mathcal{Q}_4),
\end{equation}
where
\begin{equation}
 \mathcal{Q}_c(f)(v)=\int _{\mathbb{R}^{d_v}} \int_{\mathbb{S}^{d_v-1}}B(v-v_*,\omega)[f'f_*'-ff_*]d\omega dv
\end{equation}
is the classical collision operator. The cubic terms $\mathcal{Q}_1$ -- $\mathcal{Q}_4$ are
\begin{eqnarray}
\label{Qquantum}
\left\{\begin{array}{l} 
\displaystyle\mathcal{Q}_1(f)(v)=\int _{\mathbb{R}^{d_v}} \int_{\mathbb{S}^{d_v-1}}B(v-v_*,\omega)f'f_*'f_*d\omega dv, \\
\\
\displaystyle\mathcal{Q}_2(f)(v)=\int _{\mathbb{R}^{d_v}} \int_{\mathbb{S}^{d_v-1}}B(v-v_*,\omega)f'f_*'fd\omega dv,  \\
\\
\displaystyle\mathcal{Q}_3(f)(v)=\int _{\mathbb{R}^{d_v}} \int_{\mathbb{S}^{d_v-1}}B(v-v_*,\omega)ff_*f'd\omega dv ,  \\
\\ 
\displaystyle\mathcal{Q}_4(f)(v)=\int _{\mathbb{R}^{d_v}} \int_{\mathbb{S}^{d_v-1}}B(v-v_*,\omega)ff_*f_*'d\omega dv.
\end{array}\right.
\end{eqnarray}

In order to perform the Fourier transform, we periodize the function $f$ on the domain $\mathcal{D}_L=[-L,L]^{d_v}$ ($L$ is chosen such that $L\geq \frac{3+\sqrt{2}}{2}R$, $R$ is the truncation of the collision integral which satisfies $R=2S$, where ${\mathcal B}(0,S)$ is an approximation of the support of $f$ \cite{PR00}). Using the Carleman representation \cite{Carleman33}, one can rewrite the operators as (for simplicity we only consider the 2-D Maxwellian molecules),
\begin{eqnarray}  
 \mathcal{Q}_c(f)(v)=\int _{\mathcal{B}_R} \int_{\mathcal{B}_R}\delta(x\cdot y)[f(v+x)f(v+y)-f(v+x+y)f(v)]dx dy \label{Qc}
\end{eqnarray}
and
\begin{eqnarray}  
\left\{
\begin{array}{l}
\displaystyle\mathcal{Q}_1(f)(v)=\int _{\mathcal{B}_R} \int_{\mathcal{B}_R}\delta(x\cdot y)f(v+x)f(v+y)f(v+x+y)dx dy, \\
\\
\displaystyle\mathcal{Q}_2(f)(v)=\int _{\mathcal{B}_R} \int_{\mathcal{B}_R}\delta(x\cdot y)f(v+x)f(v+y)f(v)dx dy, \\
\\
\displaystyle\mathcal{Q}_3(f)(v)=\int _{\mathcal{B}_R} \int_{\mathcal{B}_R}\delta(x\cdot y)f(v+x)f(v+x+y)f(v)dx dy, \\
\\
\displaystyle\mathcal{Q}_4(f)(v)=\int _{\mathcal{B}_R} \int_{\mathcal{B}_R}\delta(x\cdot y)f(v+y)f(v+x+y)f(v)dx dy.  
\end{array}\right.
\label{Q4}
\end{eqnarray}

Now we approximate $f$ by a truncated Fourier series,
\begin{equation}
f(v)\approx \sum_{k=-\frac{N}{2}}^{\frac{N}{2}-1}\hat{f}_ke^{i\frac{\pi}{L}k\cdot v},  \    \  \  \hat{f}_k=\frac{1}{(2L)^{d_v}}\int_{\mathcal{D}_L}f(v)e^{-i\frac{\pi}{L}k\cdot v}dv.
\end{equation}
Plugging it into (\ref{Qc}) (\ref{Q4}), one can get the $k$-th mode of $\hat{\mathcal{Q}}_q$. The classical part is the same as those in the previous method \cite{MP06}. We will mainly focus on the cubic terms. 

Define the kernel modes
\begin{equation}
\beta(l,m)=\int_{\mathcal{B}_R}\int_{\mathcal{B}_R} \delta(x\cdot y)e^{i\frac{\pi}{L}l\cdot x}e^{i\frac{\pi}{L}m\cdot y}dxdy.
\end{equation}
Following \cite{MP06}, $\beta(l,m)$ can be decomposed as
\begin{equation}
\beta(l,m)=\frac{\pi}{M}\sum_{p=0}^{M-1}\alpha_p(l)\alpha'_p(m)
\end{equation}
with
\begin{equation}
\alpha_p(l)=\phi(l\cdot (\cos \theta_p,\sin \theta_p)),  \   \    \   \alpha_p'(m)=\phi(m\cdot (-\sin \theta_p,\cos \theta_p)),
\end{equation}
where $\phi(s)=\frac{2L}{\pi s}\sin(\frac{\pi}{L}Rs)$, $M$ is the number of equally spaced points in $[0, \frac{\pi}{2}]$ and $\theta_p=\frac{\pi}{2}\frac{p}{M}$. Then

\begin{itemize}
\item The $k$-th coefficient of $\hat{\mathcal{Q}}_1$ is
\begin{eqnarray}
\sum_{\substack{l,m,n=-\frac{N}{2}\\l+m+n=k}}^{\frac{N}{2}-1}\beta(l+n,m+n)\hat{f}_l\hat{f}_m\hat{f}_n&=&\frac{\pi}{M}\sum_{p=0}^{M-1}\sum_{n=-\frac{N}{2}}^{\frac{N}{2}-1}\left[\sum_{\substack{l,m=-\frac{N}{2}\\ l+m=k-n}}^{\frac{N}{2}-1}\alpha_p(l+n)\alpha'_p(m+n)\hat{f}_l\hat{f}_m\right]\hat{f}_n \nonumber \\
&=&\frac{\pi}{M}\sum_{p=0}^{M-1}\sum_{n=-\frac{N}{2}}^{\frac{N}{2}-1} \hat{g}_{k-n}(n)\hat{f}_n.
\end{eqnarray}
Terms inside the bracket is a convolution (defined as $\hat{g}_{k-n}(n)$), which can be computed by the Fast Fourier Transform (FFT). However, the outside structure is not a convolution, since $\hat{g}_{k-n}(n)$ itself depends on $n$. So we compute this part directly.

\item The $k$-th coefficient of $\hat{\mathcal{Q}}_2$ is
\begin{eqnarray}
\sum_{\substack{l,m,n=-\frac{N}{2}\\l+m+n=k}}^{\frac{N}{2}-1}\beta(l,m)\hat{f}_l\hat{f}_m\hat{f}_n=\frac{\pi}{M}\sum_{p=0}^{M-1}\sum_{n=-\frac{N}{2}}^{\frac{N}{2}-1}\left[\sum_{\substack{l,m=-\frac{N}{2}\\ l+m=k-n}}^{\frac{N}{2}-1}\alpha_p(l)\alpha'_p(m)\hat{f}_l\hat{f}_m\right]\hat{f}_n .
\end{eqnarray}
In this case, both inside and outside are convolutions. The FFT can be implemented easily.

\item The $k$-th coefficient of $\hat{\mathcal{Q}}_3$ is
\begin{eqnarray}
\sum_{\substack{l,m,n=-\frac{N}{2}\\l+m+n=k}}^{\frac{N}{2}-1}\beta(l+m,m)\hat{f}_l\hat{f}_m\hat{f}_n=\frac{\pi}{M}\sum_{p=0}^{M-1}\sum_{n=-\frac{N}{2}}^{\frac{N}{2}-1}\alpha_p(l+m)\left[\sum_{\substack{l,m=-\frac{N}{2}\\ l+m=k-n}}^{\frac{N}{2}-1}\alpha'_p(m)\hat{f}_l\hat{f}_m\right]\hat{f}_n.
\end{eqnarray}
Factoring out $\alpha_p(l+m)$, both inside and outside are convolutions again.

\item The $k$-th coefficient of $\hat{\mathcal{Q}}_4$ is
\begin{eqnarray}
\sum_{\substack{l,m,n=-\frac{N}{2}\\l+m+n=k}}^{\frac{N}{2}-1}\beta(m,l+m)\hat{f}_l\hat{f}_m\hat{f}_n=\frac{\pi}{M}\sum_{p=0}^{M-1}\sum_{n=-\frac{N}{2}}^{\frac{N}{2}-1}\alpha_p'(l+m)\left[\sum_{\substack{l,m=-\frac{N}{2}\\ l+m=k-n}}^{\frac{N}{2}-1}\alpha_p(m)\hat{f}_l\hat{f}_m\right]\hat{f}_n .
\end{eqnarray}
This term can be evaluated similarly as $\hat{\mathcal{Q}}_3$.
\end{itemize}

\begin{remark}
The computational cost of this quantum solver is $O(MN^4\log N)$, which mainly comes from computing $\mathcal{Q}_1$. This cost is higher than $O(MN^4)$ of the discrete velocity model. But taking into account the high accuracy and small value of $\log N$ ($N$ is not very big in the real simulation), our method is still 
more attractive than the quadrature method. The fast algorithm for the quantum collision operator remains an open problem.
\end{remark}

\subsection{Numerical Accuracy}
To illustrate the accuracy of the above method, we test it on a steady state, namely, we compute $\mathcal{Q}_q(\mathcal{M}_q)$ and check its max norm. In all the numerical simulations, the particles are assumed to be the 2-D Maxwellian molecules.

\begin{figure}[h!]
\centering
\includegraphics[width=15cm]{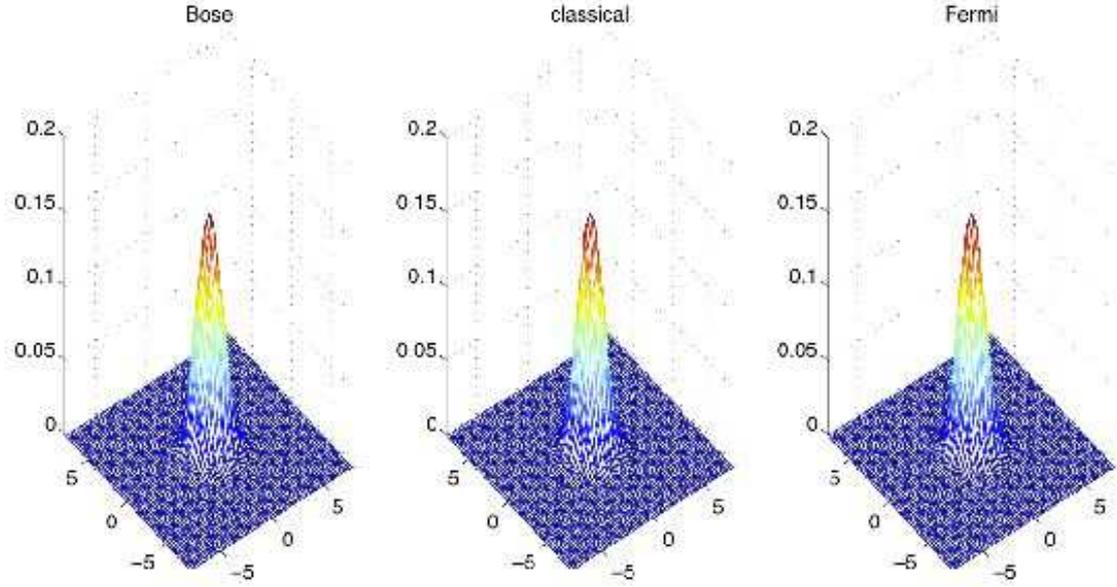}
\caption{The Maxwellians at $\rho=1$, $T=1$, $\theta_0=0.01$. Left: Bose gas; Center: classical gas; Right: Fermi gas.
}
\label{com_s}
\end{figure}

\begin{figure}[h!]
\centering
\includegraphics[width=15.cm]{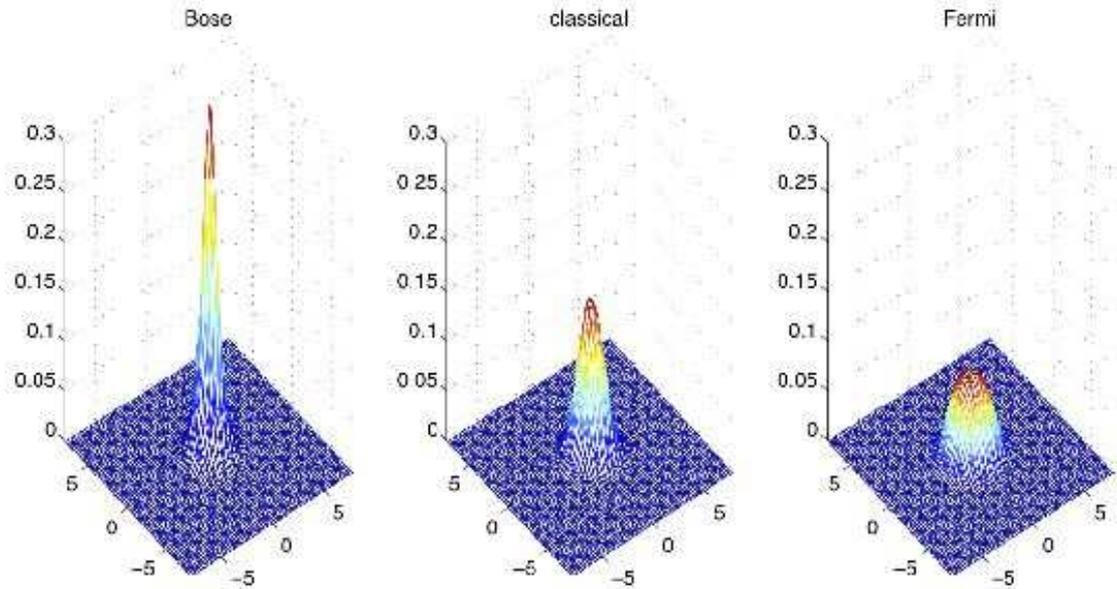}
\caption{The Maxwellians at $\rho=1$, $T=1$, $\theta_0=9$. Left: Bose gas; Center: classical gas (same as in Fig.\ref{com_s}); Right: Fermi gas.
}
\label{com_l}
\end{figure}

Let $\rho=1$, $T=1$, from (\ref{rhoue}) one can adjust $\theta_0$ to get $z$ that lies in different physical regimes. When $\theta_0=0.01$ ($\hbar=0.1$), $z_{\text{Bose}}=0.001590$, $z_{\text{Fermi}}=0.001593$. In this situation, the quantum effect is very small. The Maxwellians for the Bose gas, classical gas and Fermi gas are almost the same (Fig.\ref{com_s}). When we increase $\theta_0$, say $\theta_0=9$ ($\hbar=3$), $z_{\text{Bose}}=0.761263$, $z_{\text{Fermi}}=3.188717$, the difference between the quantum gases and the classical gas is evident (Fig.\ref{com_l}).

In Table \ref{table1}, we list the values of $\parallel\mathcal{Q}_c(\mathcal{M}_c)\parallel_{L^\infty}$ and $\parallel\mathcal{Q}_q(\mathcal{M}_q)\parallel_{L^{\infty}}$ computed on different meshes N=16, 32, 64 (number of points in $v$ direction), M=4 (number of points in angular direction $\theta_p$; it is not necessary to put too many points since $M$ won't effect the spectral accuracy, see \cite{MP06}). The computational domain is $[-8,8]\times[-8,8]$ ($L=8$).

\begin{table}[h!]
\centering
 \begin{tabular}{c  l | c | c | c | c}
\hline 
\\
&  & $16\times 16$ & $32\times 32$ & $64\times 64$  & convergence rate\\
 \hline
\\
 classical gas   &   &2.1746e-04& 3.8063e-12& 1.9095e-16 &20.0253 \\
 \hline
\\
 Bose gas  &$\theta_0=0.01$& 2.1084e-04& 2.5512e-10& 1.9080e-16 & 20.0036\\
                &$\theta_0=9$  & 0.4891 & 0.0310 & 1.3496e-04 &5.9117\\
                &  $\theta_0=9$, $L=6$                        & 0.1815 & 0.0052 &  4.0278e-06 & 7.7298\\
 \hline
\\
 Fermi  gas &$\theta_0=0.01$ & 2.2397e-04 & 1.6485e-10& 1.9152e-16 & 20.0445\\
                     & $\theta_0=9$& 8.9338e-04&2.0192e-06&1.5962e-10 & 11.2081\\
 \hline
   \end{tabular}
\caption{Comparison of the quantum collision solver on different Maxwellians ($L=8$ unless specified).
}
\label{table1}
\end{table}

These results confirm the spectral accuracy of the method, although the accuracy in the quantum regime is not as good as that in the classical regime. This is because the regularity of the quantum Maxwellians becomes worse when $\theta_0$ is increasing, or strictly speaking, the mesh size $\Delta v$ is not small enough to capture the shape of the Maxwellians. To remedy this problem, one can add more grid points or more effectively, shorten the computational domain. For the Bose-Einstein distribution, we also include the results computed on $[-6,6]\times[-6,6]$ in Table \ref{table1}. One can clearly see the improvements.

\subsection{Relaxation to Equilibrium}
Let us consider the space homogeneous quantum Boltzmann equation for the 
2-D Maxwellian molecules. As already mentioned, this equation satisfies the entropy condition,
 and the equilibrium states are the entropy minimizers.  Hence, we first consider the quantum Boltzmann equation for a Fermi gas with an initial datum $0\leq f_0 \leq \frac{1}{\theta_0}$ and observe the relaxation to equilibrium of the distribution function. Then, we take a Bose gas for which the entropy is now sublinear and fails to prevent concentration, which is consistent with the fact that condensation may occur in the long-time limit.

\paragraph{Fermi gas.} The initial data is chosen as the sum of two Maxwellian functions
\begin{equation}
\label{my:data}
f_0(v) =   \exp\left(-\frac{|v-v_1|^2}{2}\right) +  \exp\left(-\frac{|v+v_1|^2}{2}\right);\quad v\in \mathbb{R}^{2},
\end{equation}
with $v_1 = (2,1)$. The final time of the simulation is $T_{end} = 0.5$, which is very
close to the stationary state. 

In the spatially homogeneous setting, Pauli's exclusion principle facilitates things  because of the additional $L^\infty$ bound $0\leq f(t)\leq \frac{1}{\theta_0}$. In this case, the convergence to equilibrium in a  weak sense has been shown by Lu \cite{Lu01}. Later Lu and Wennberg proved the strong $L^1$ stability \cite{LW03}.  However, no constructive result in this direction has ever been obtained, neither has any entropy-dissipation inequality been established. 

\begin{figure}[h!]
\centering
\includegraphics[width=15.cm]{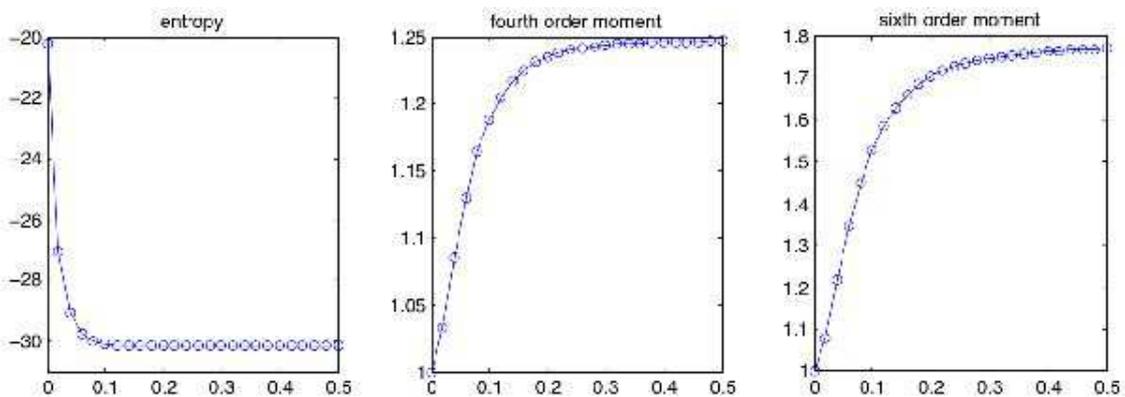}
\caption{Fermi gas. Time evolution of the entropy, fourth and sixth order moments.}
\label{ferm1d}
\end{figure}

In Fig.\ref{ferm1d} we report the time evolution of the entropy and  the fourth and sixth order moments of the distribution with respect to the velocity variable. We indeed observe the convergence to a steady state of the entropy and also of high order moments  when $t\rightarrow \infty$.

\begin{figure} [h!]
\centering
\includegraphics[width=7.25cm]{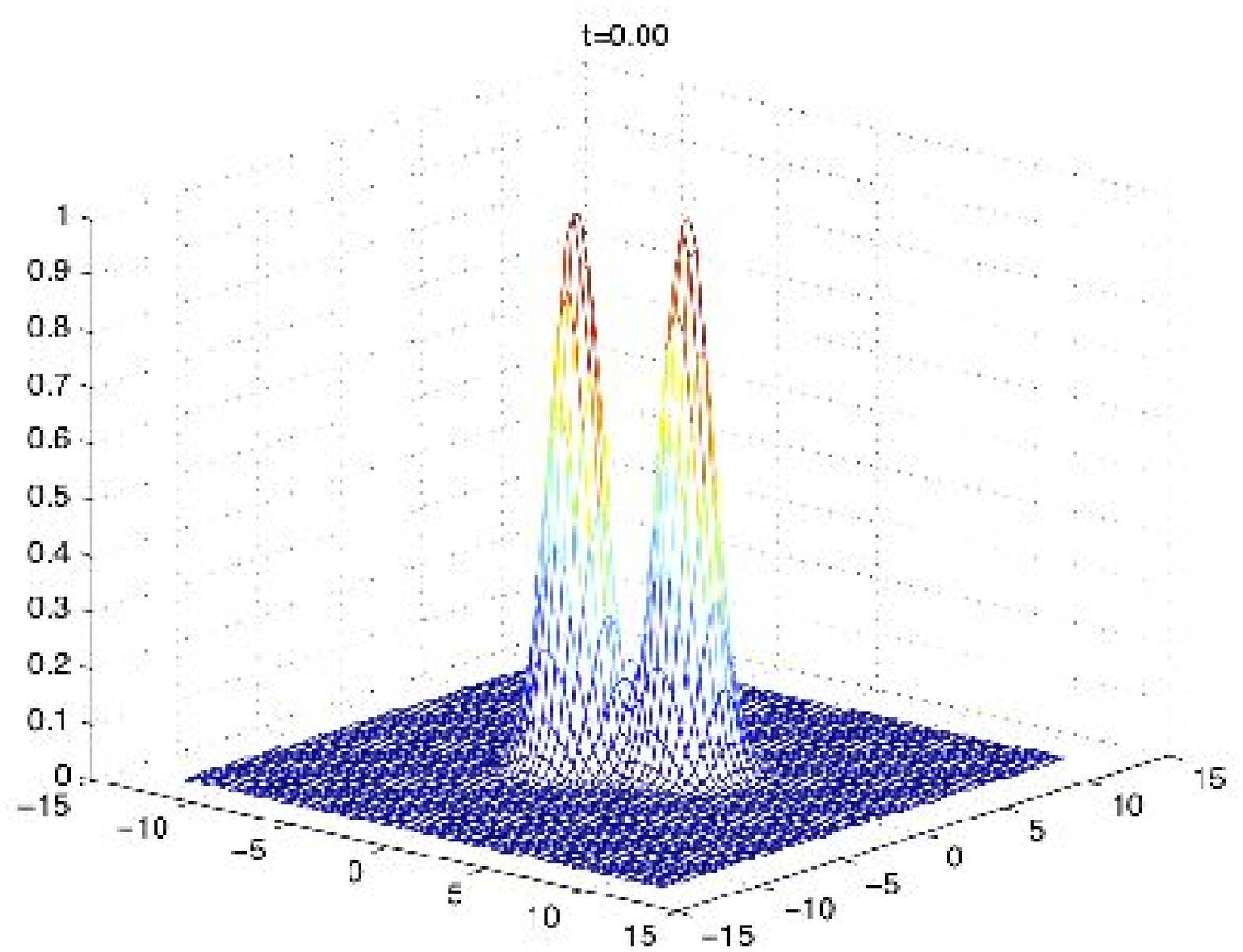}
\includegraphics[width=7.25cm]{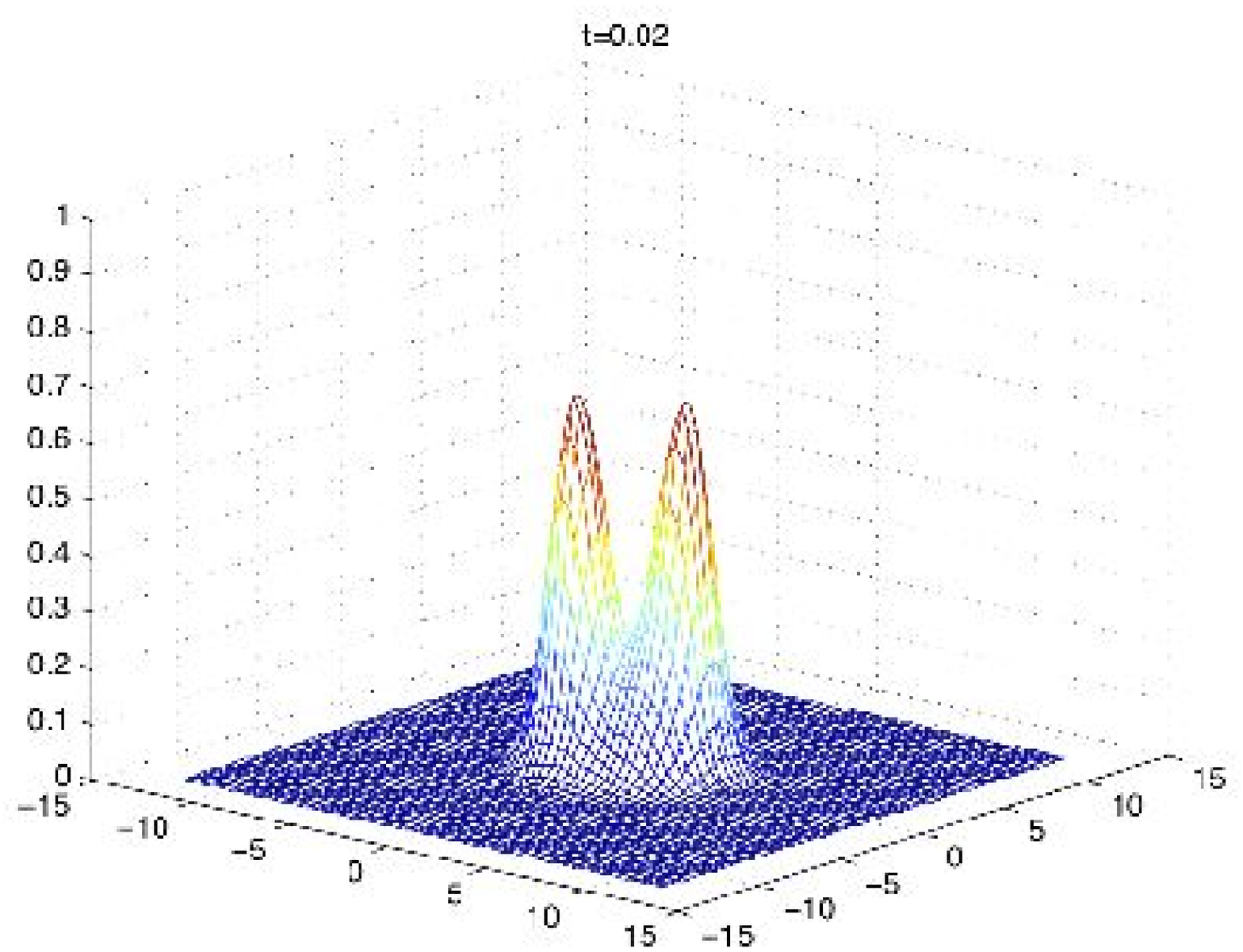}
\includegraphics[width=7.25cm]{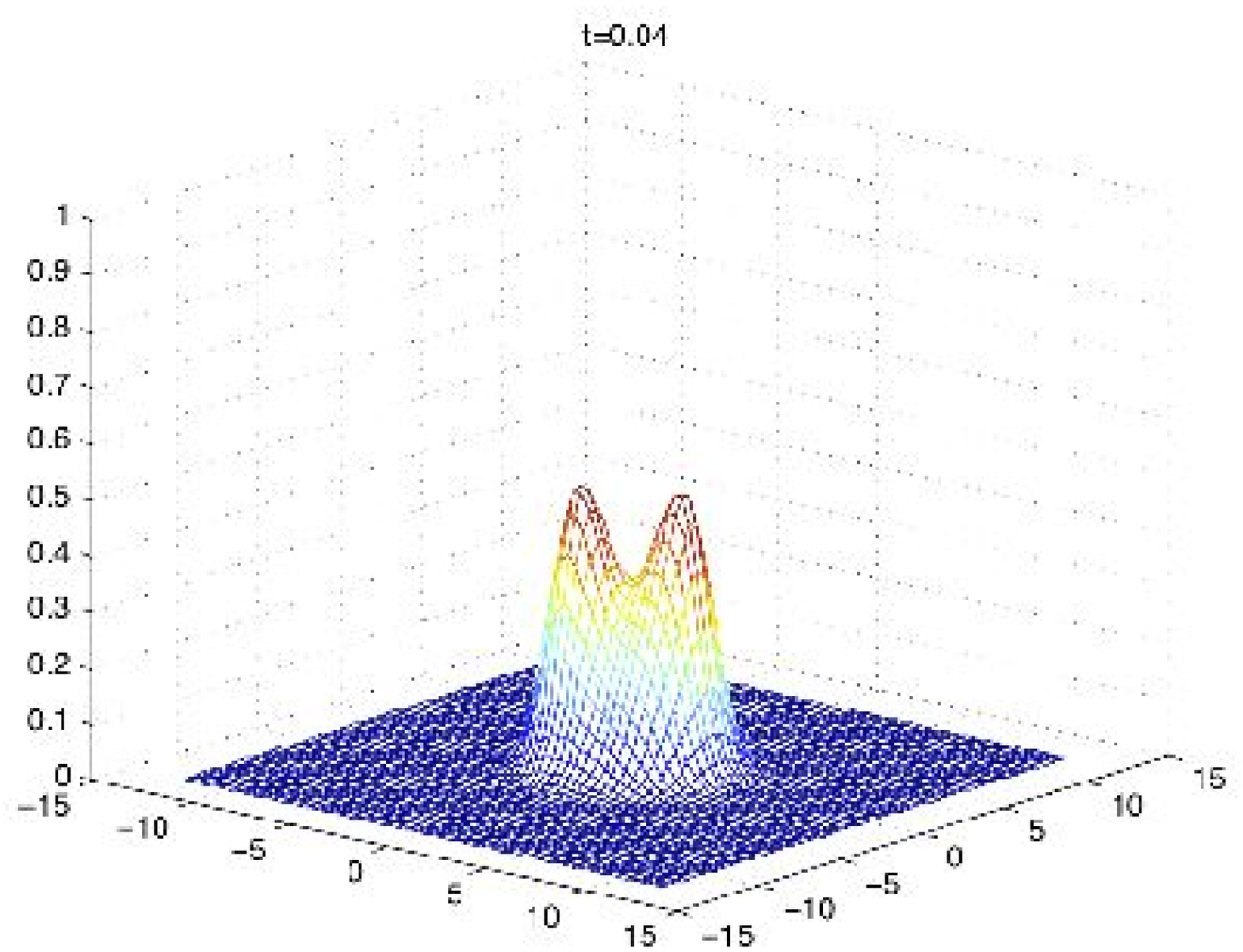}
\includegraphics[width=7.25cm]{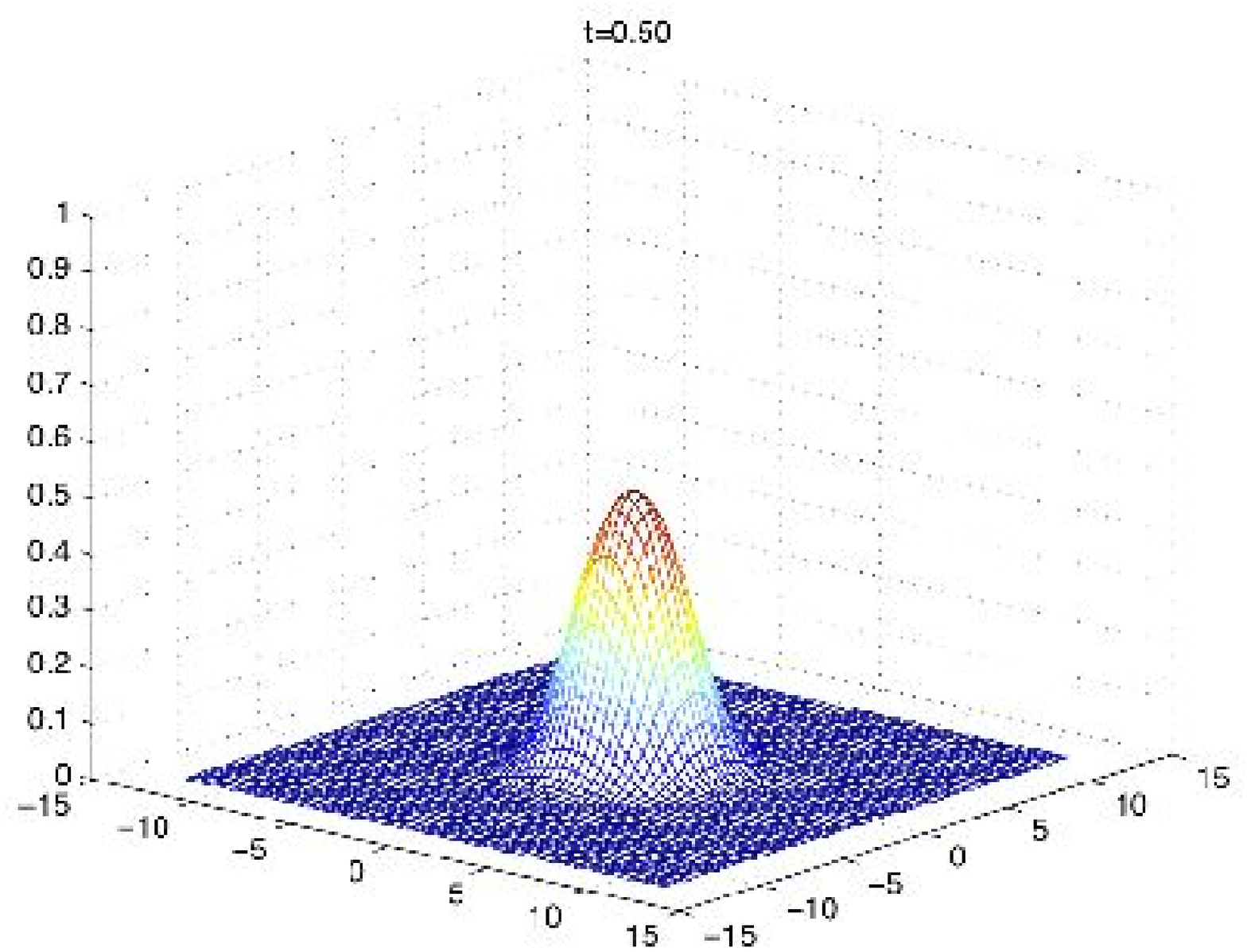}
\caption{Fermi gas. Time evolution of the distribution function $f(t,v_x,v_y)$ 
with $N=64$ modes  at times $t=0, 0.02, 0.04$ and $0.5$.}
\label{fv2D}
\end{figure}

In  Fig.\ref{fv2D} we also report the time evolution of the level set of the distribution function $f(t,v_x,v_y)$ obtained  with $N=64$ modes at different times. Initially the level set of the initial data corresponds to two spheres  in the velocity space. Then, the two distributions start to mix together  until the stationary state is reached, represented by a single  centered sphere. It is clear that the spherical shapes of the level  sets are described with great accuracy by the spectral method.

\paragraph{Bose gas.} This is an even more challenging problem since there is no convergence result, due to the lack of a priori bound.  Lu \cite{Lu00} has attacked this problem with the well-developed tools of the modern spatially homogeneous theory and proved that the solution (with a very low temperature) converges to equilibrium in a weak sense. 
In \cite{EM01}, the authors studied  an one dimensional model and proved existence theorems, and convergence to a Bose distribution having a singularity when time goes to infinity because Bose condensation cannot occur in finite time.

Here we investigate the convergence to equilibrium for space homogeneous model in 2-D, for which condensation cannot occur. We consider the following initial datum 
\begin{equation}
f_0(v) =   \frac{1}{4\pi\,T_0}\exp\left(-\frac{|v-v_1|^2}{2T_0}\right) +  \exp\left(-\frac{|v+v_1|^2}{2T_0}\right);\quad v\in \mathbb{R}^{2},
\end{equation}
with $v_1 = (1,1/2)$ and $T_0=1/4$. 

\begin{figure}[h!]
\centering
\includegraphics[width=15cm]{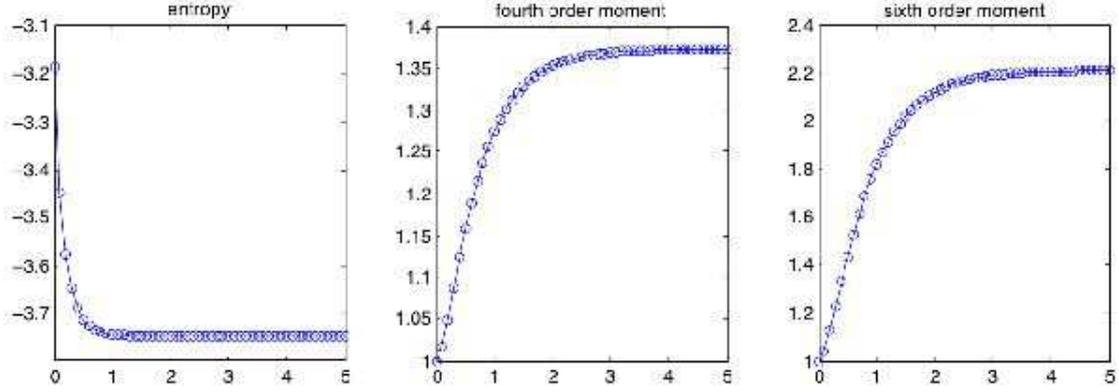}
\caption{Bose gas. Time evolution of the entropy, fourth and sixth order moments.}
\label{bose1d}
\end{figure}

We still observe the convergence to equilibrium and convergence of high order moments when $t \rightarrow \infty$ in Fig.\ref{bose1d}.

\begin{figure} [h!]
\centering
\includegraphics[width=7.25cm]{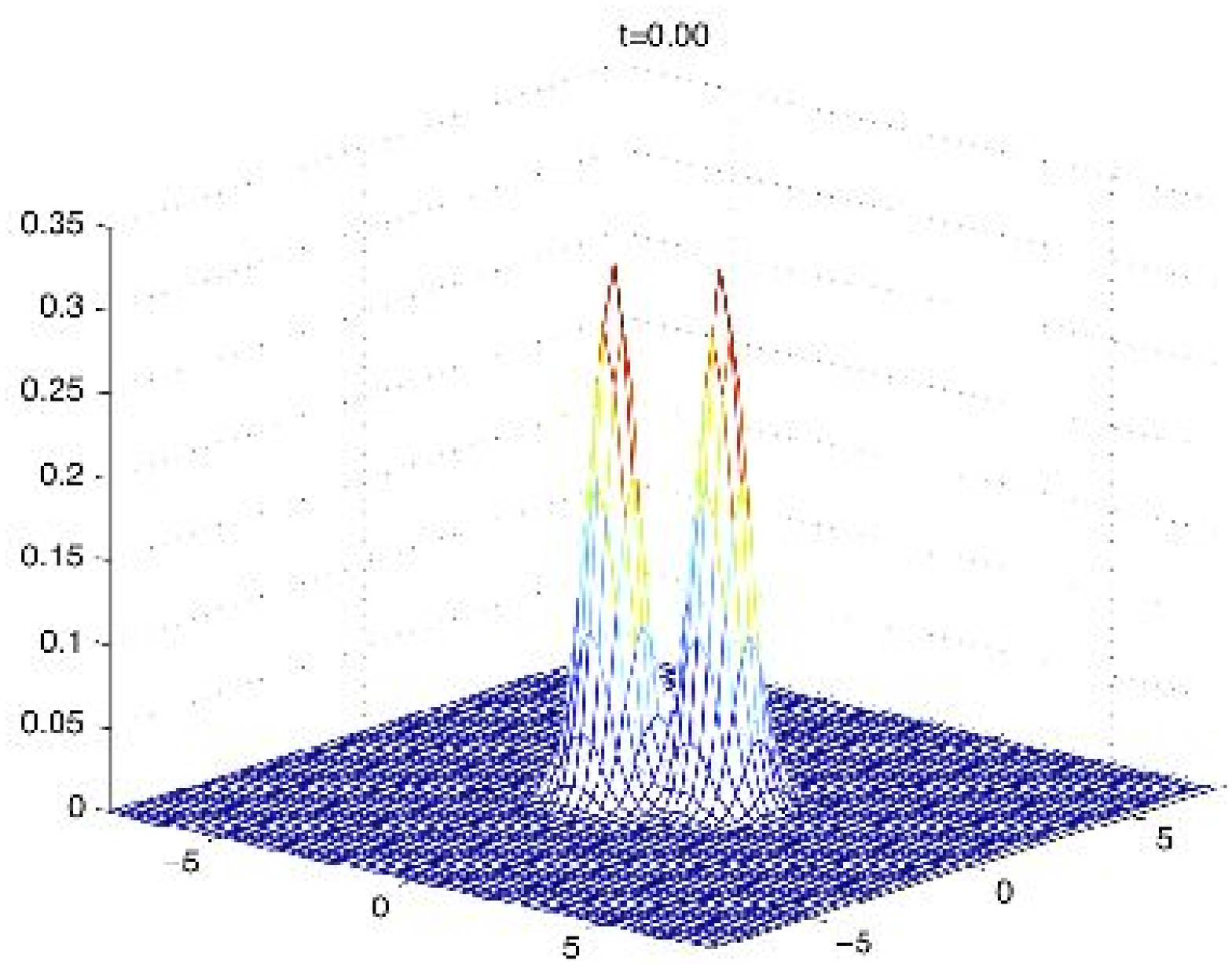}
\includegraphics[width=7.25cm]{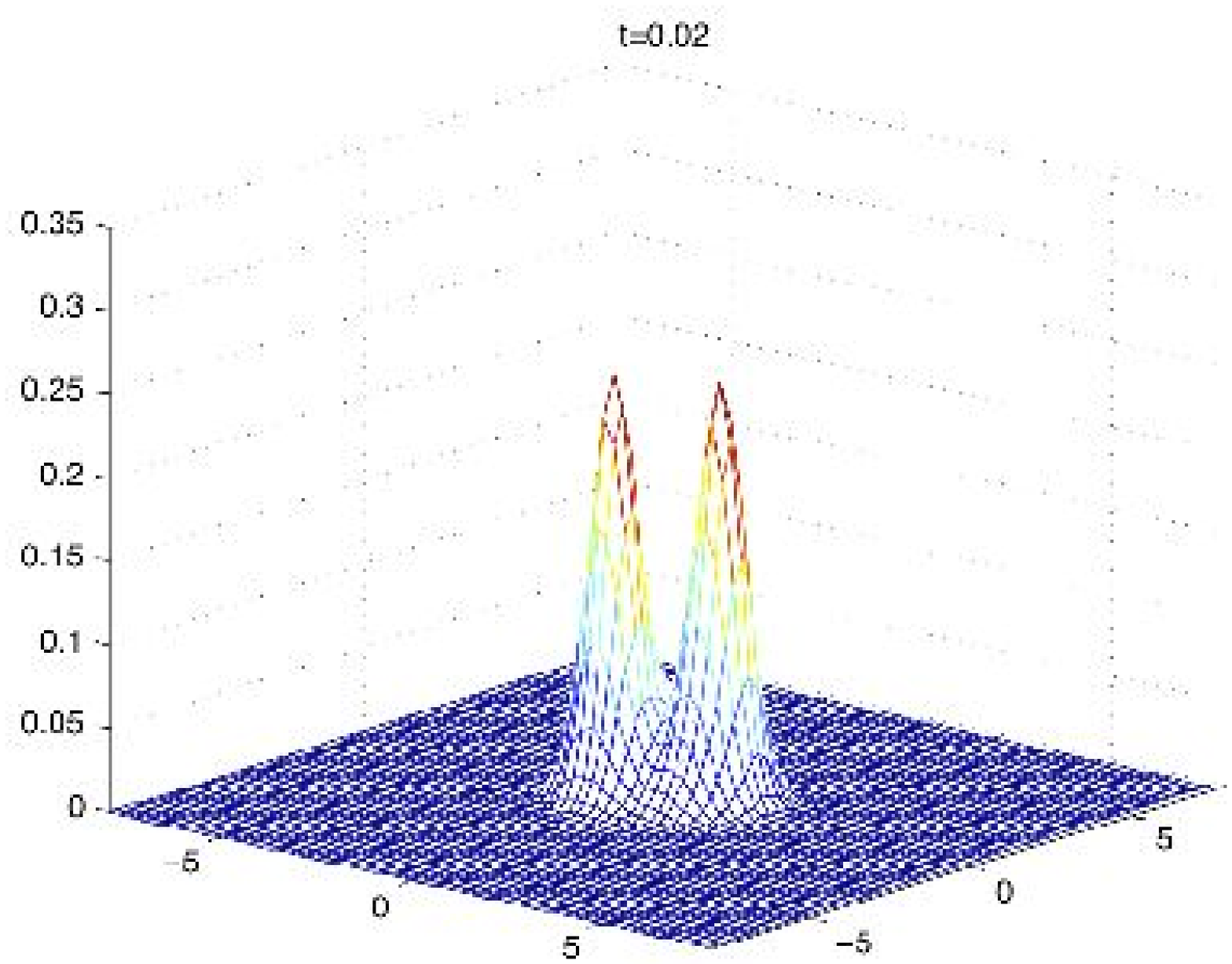}
\includegraphics[width=7.25cm]{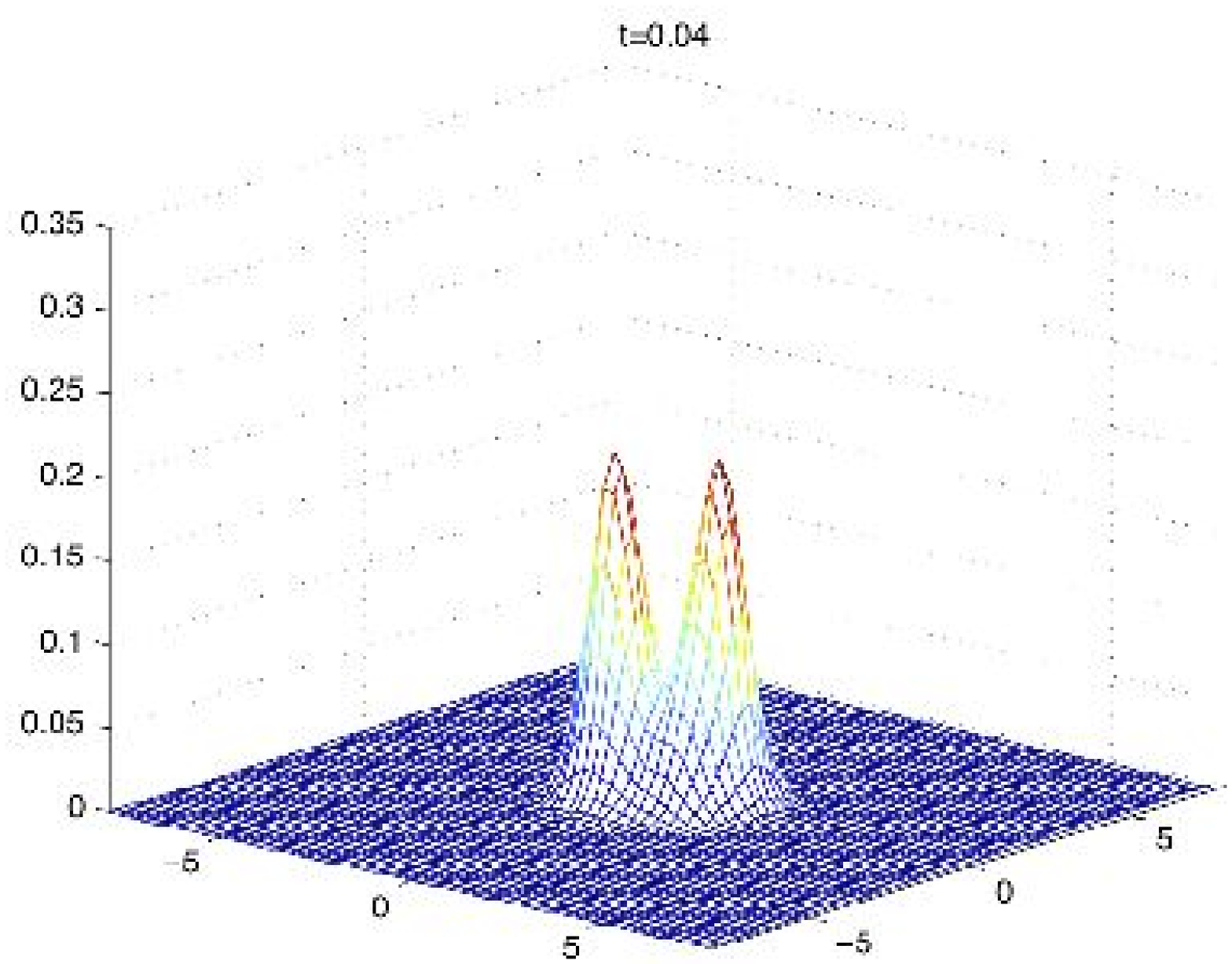}
\includegraphics[width=7.25cm]{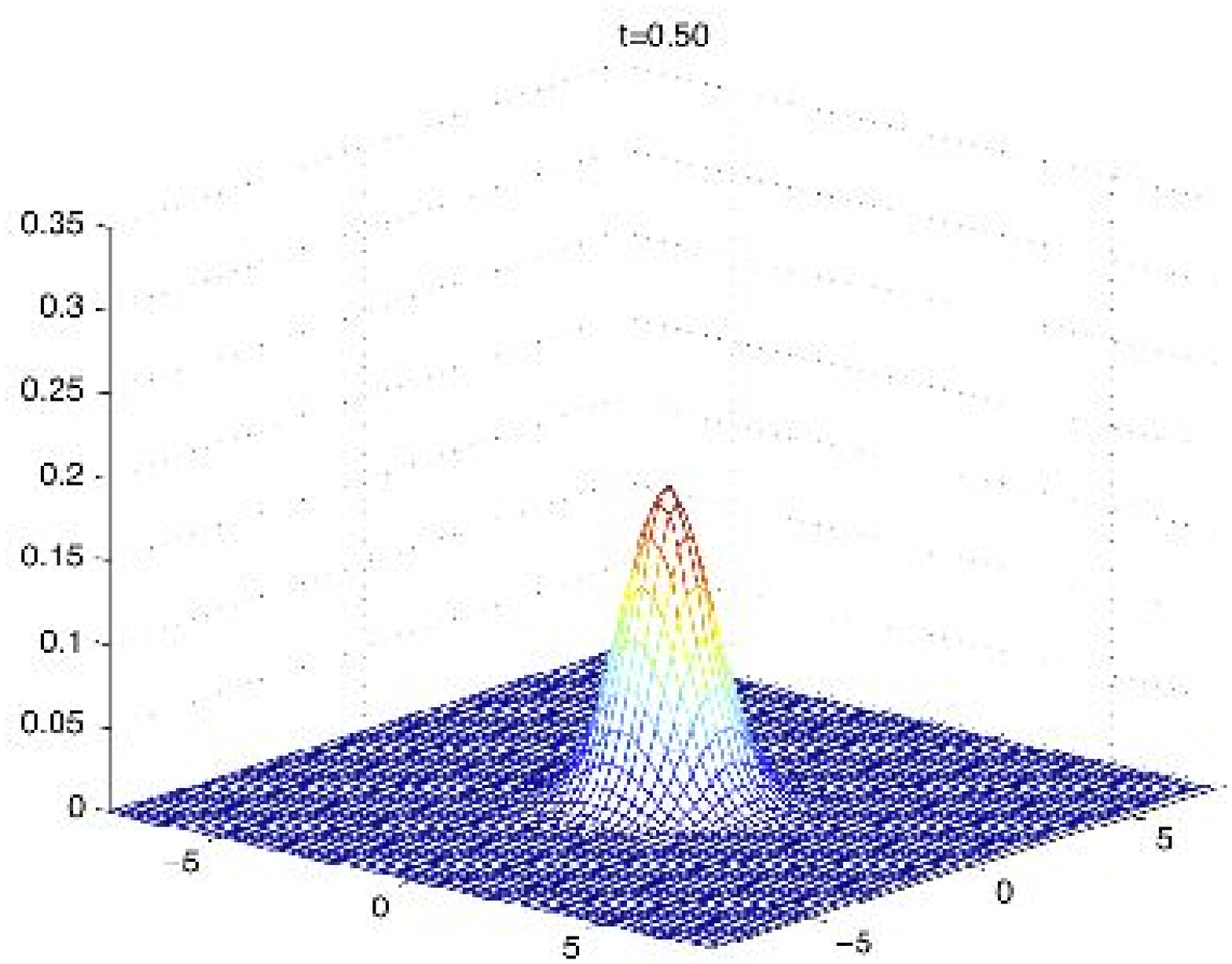}
\caption{Bose gas. Time evolution of the distribution function $f(t,v_x,v_y)$ 
with $N=64$ modes  at times $t=0, 0.02, 0.04$ and $0.5$.}
\label{bv2D}
\end{figure}

In  Fig.\ref{bv2D} we report the time evolution of the level set of the distribution function $f(t,v_x,v_y)$ obtained  with $N=64$ modes at different times and observe the trend to equilibrium.

\section{A Scheme Efficient in the Fluid Regime}
So far we have only considered spatially homogeneous quantum Boltzmann equations, now what happens for spatially inhomogeneous data? Due to the natural bound $0\leq f(t) \leq \frac{1}{\theta_0}$, the Boltzmann-Fermi model seems to be well understood mathematically \cite{Villani02}. The situation is completely different for the Boltzmann-Bose model, since  singular measures may occur \cite{Villani02}.

We first review the scheme in \cite{FJ10} for the classical Boltzmann equation
\begin{equation} 
\frac{\partial f}{\partial t}\,+\,v\cdot \nabla _x f=\frac{1}{\epsilon}\mathcal{Q}_c(f).
\end{equation}
The first-order scheme reads:
\begin{equation}\label{scheme}
\frac{f^{n+1}-f^n}{\Delta t}+v\cdot \nabla_xf^n=\frac{\mathcal{Q}_c(f^n)-\lambda(\mathcal{M}_c^n-f^n)}{\epsilon}+\frac{\lambda(\mathcal{M}_c^{n+1}-f^{n+1})}{\epsilon},\end{equation}
where $\lambda$ is some appropriate approximation of $|\nabla \mathcal{Q}_c|$ (can be made time dependent). To solve $f^{n+1}$ explicitly, we need to compute $\mathcal{M}_c^{n+1}$ first. Since the right hand side of (\ref{scheme}) is conservative, it vanishes when we take the moments (multiply by $\phi(v)=(1,v,\frac{1}{2}v^2)^T$ and integrate with respect to $v$). Then (\ref{scheme}) becomes
\begin{equation}
\frac{U^{n+1}-U^n}{\Delta t}+\int \phi(v) v\cdot \nabla_x f^ndv=0,
\end{equation}
where $U=(\rho,\rho u, \rho e+\frac{1}{2}\rho u^2)^T$ is the conserved quantities. Once we get $U^{n+1}$, $\mathcal{M}_c^{n+1}$ is known. Now $f^{n+1}$ in (\ref{scheme}) is easy to obtain.

When generalizing the above idea to the quantum Boltzmann equation (\ref{QBE}), the natural idea is to replace $\mathcal{Q}_c$ and $\mathcal{M}_c$ in (\ref{scheme}) by $\mathcal{Q}_q$ and $\mathcal{M}_q$ respectively. However, as mentioned in section 2, one has to invert the nonlinear system (\ref{rhoue})  to get $z$ and $T$. Experiments show that the iterative methods do converge when the initial guess is close to the solution (analytically, this system has a solution \cite{AL97}). But how to set a good initial guess for every spatial point and every time step is not an easy task, especially when $\rho$ and $e$ are not continuous.

Here we propose to use a `classical' BGK operator to penalize $\mathcal{Q}_q$. Specifically, we replace the temperature $T$ with the internal energy $e$ in the classical Maxwellian using relation $e=\frac{d_v}{2}T$ (true for classical monatomic gases) and get
\begin{equation}
\mathcal{M}_c=\frac{\rho}{(2\pi T)^{\frac{d_v}{2}}}e^{-\frac{(v-u)^2}{2T}}=\rho \left(\frac{d_v}{4\pi e} \right)^{\frac{d_v}{2}}e^{-\frac{d_v}{4e}(v-u)^2}.
\end{equation}
An important property of $\mathcal{M}_c$ is that it has the same first five moments as $\mathcal{M}_q$.

Now our new scheme for QBE (\ref{QBE}) can be written as
\begin{equation} \label{scheme1}
\frac{f^{n+1}-f^n}{\Delta t}+v\cdot \nabla_xf^n=\frac{\mathcal{Q}_q(f^n)-\lambda(\mathcal{M}_c^n-f^n)}{\epsilon}+\frac{\lambda(\mathcal{M}_c^{n+1}-f^{n+1})}{\epsilon}.\end{equation}
Since the right hand side is still conservative, one computes $\mathcal{M}_c^{n+1}$ the same as for (\ref{scheme}).

It is important to notice that $z$ and $T$ are not present at all in this new scheme, thus one does not need to invert the 2 by 2 system (\ref{rhoue})  during the time evolution. If they are desired variables for output, one only needs to convert between $\rho$, $e$ and $z$, $T$ at the final output time.

\subsection{Asymptotic Property of the New Scheme}
In this subsection we show that the new scheme, when applied to the quantum BGK equation, has the property (\ref{property}). Consider the following time discretization:
\begin{equation} \label{BGK}
\frac{f^{n+1}-f^n}{\Delta t}+v\cdot \nabla_xf^n=\frac{(\mathcal{M}_q^n-f^n)-\lambda(\mathcal{M}_c^n-f^n)}{\epsilon}+\frac{\lambda(\mathcal{M}_c^{n+1}-f^{n+1})}{\epsilon}.\end{equation}
Some simple mathematical manipulation on (\ref{BGK}) gives
\begin{equation} \label{eps}
f^{n+1}-\mathcal{M}_q^{n+1}=\frac{1+(\lambda-1)\frac{\Delta t}{\epsilon}}{1+\lambda\frac{\Delta t}{\epsilon}}(f^n-\mathcal{M}_q^n)-\frac{\Delta t}{1+\lambda\frac{\Delta t}{\epsilon}}v\cdot \nabla_xf^n+(\mathcal{M}_q^n-\mathcal{M}_q^{n+1})+\frac{\lambda \frac{\Delta t}{\epsilon}}{1+\lambda \frac{\Delta t}{\epsilon}}(\mathcal{M}_c^{n+1}-\mathcal{M}_c^n).
\end{equation}
Assume all the functions are smooth. When $\lambda >\frac{1}{2}$,
\begin{equation}
|f^{n+1}-\mathcal{M}_q^{n+1}|\leq \alpha|f^n-\mathcal{M}_q^n|+O(\epsilon+\Delta t),
\end{equation}
where $0<\alpha=|1+(\lambda-1)\frac{\Delta t}{\epsilon}|/|1+\lambda\frac{\Delta t}{\epsilon}|<1$ uniformly in $\epsilon$ and $\Delta t$. The $O(\epsilon)$ term comes from the second term of the right hand side of (\ref{eps}). The $O(\Delta t)$ term is from the third and fourth terms.
Then
\begin{equation}
|f^n-\mathcal{M}_q^n|\leq \alpha^n|f^0-\mathcal{M}_q^0|+O(\epsilon+\Delta t)\,.
\end{equation}
Since $\Delta t$ is taken bigger than $\epsilon$, this implies the property (\ref{property}). It is interesting to point out that $f$ approaches $\mathcal{M}_q$, not $\mathcal{M}_c$, with (\ref{BGK}).

\begin{remark} The first order (in-time) method can be extended to a second order by an Implicit-Explicit (IMEX) method (see also \cite{FJ10}):
\begin{equation} \left\{
\begin{aligned}
&\frac{f^*-f^n}{\Delta t/2}+v\cdot \nabla_xf^n=\frac{\mathcal{Q}_q(f^n)-\lambda(\mathcal{M}_c^n-f^n)}{\epsilon}+\frac{\lambda(\mathcal{M}^*_c-f^*)}{\epsilon},\\
&\frac{f^{n+1}-f^n}{\Delta t}+v\cdot \nabla_xf^*=\frac{\mathcal{Q}_q(f^*)-\lambda(\mathcal{M}_c^*-f^*)}{\epsilon}+\frac{\lambda(\mathcal{M}^n_c-f^n)+\lambda(\mathcal{M}_c^{n+1}-f^{n+1})}{2\epsilon}.\end{aligned}
\right.
\end{equation}
This scheme can be shown to have the same property (\ref{property}) on the quantum BGK equation.
\end{remark}

\section{Numerical Examples}
In this section, we present some numerical results of our new scheme (\ref{scheme1}) (a second order finite volume method with slope limiters \cite{Leveque} is applied to the transport part) on the 1-D shock tube problem. The initial condition is

\begin{equation} \label{IC}
 \left\{ 
\begin{array}{ll}
(\rho_l,u_l,T_l)=(1,0,1) &   \mbox{if $0\leq x \leq 0.5$},\\
(\rho_r,u_r,T_r)=(0.125,0,0.25)&   \mbox{if $0.5< x \leq 1$}.
\end{array} 
\right. 
\end{equation}

The particles are again assumed to be the 2-D Maxwellian molecules and we adjust $\theta_0$ to get different initial data for both the Bose gas and the Fermi gas. 

In all the regimes, besides the directly computed macroscopic quantities, we will show the fugacity $z$ and temperature $T$ as well. They are computed as follows. First,  (\ref{rhoue}) ($d_v=2$) leads to
\begin{equation} \label{inv}
\frac{Q_1^2(z)}{Q_2(z)}=\frac{\theta_0}{2\pi}\frac{\rho}{e}.
\end{equation}
We treat the left hand side of (\ref{inv}) as one function
of $z$, and invert it by the secant method. Once  $z$ is obtained, $T$ can be 
computed easily using for example the first equation of (\ref{rhoue}). To evaluate the quantum function $Q_{\nu}(z)$, the expansion (\ref{expan1}) is used for the Bose-Einstein function. The Fermi-Dirac function is computed by a direct numerical integration. The approach adopted here is taken from \cite{PTVF} (Chapter 6.10). 

When approximating the collision operator $\mathcal{Q}_q$, we always take $M=4$, $N=32$ and $L=8$, except $L=6$ for the Bose gas in the quantum regime.

\subsection{Hydrodynamic Regime}
\begin{figure} [h!]
\centering
\includegraphics[width=7.25cm]{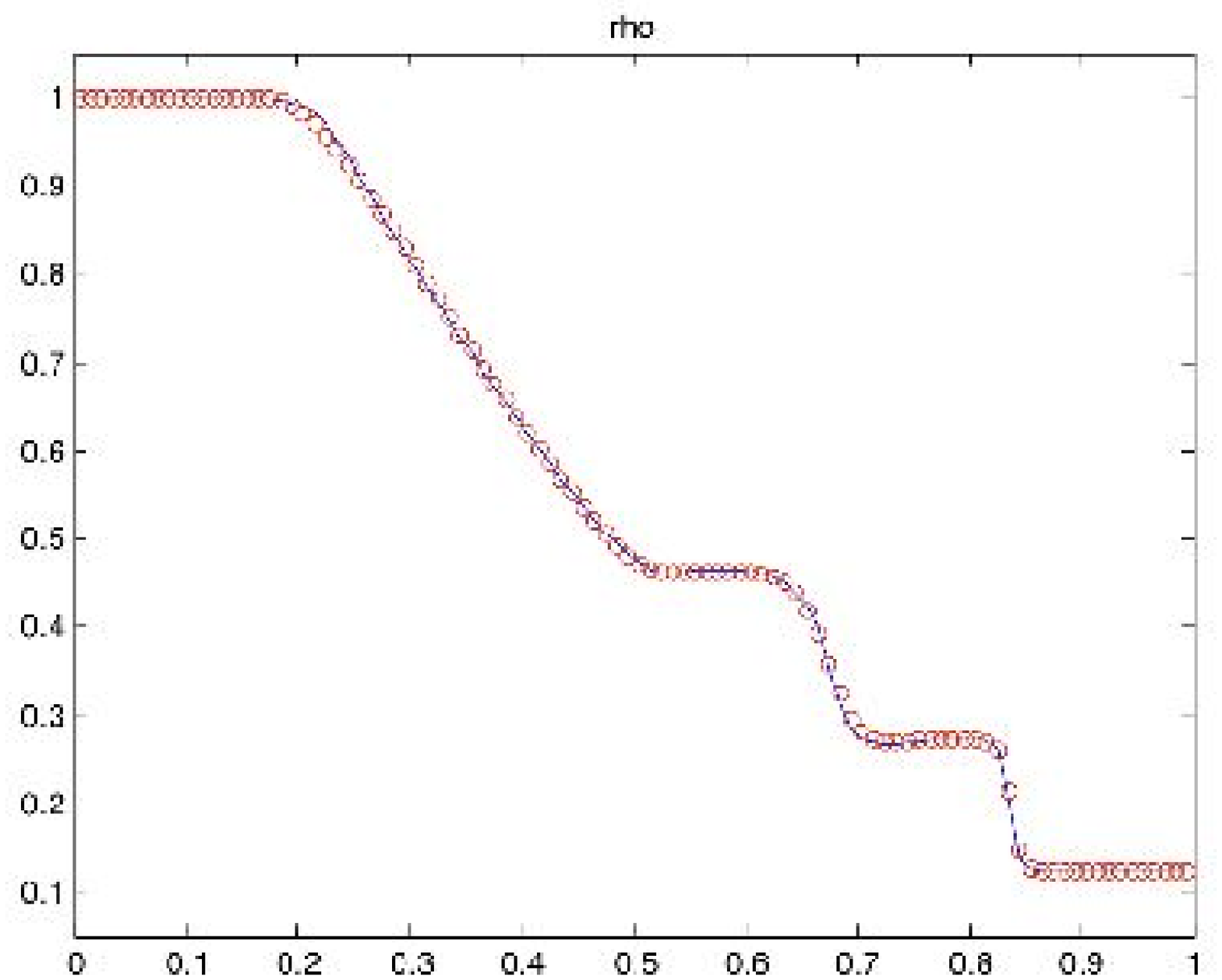}
\includegraphics[width=7.25cm]{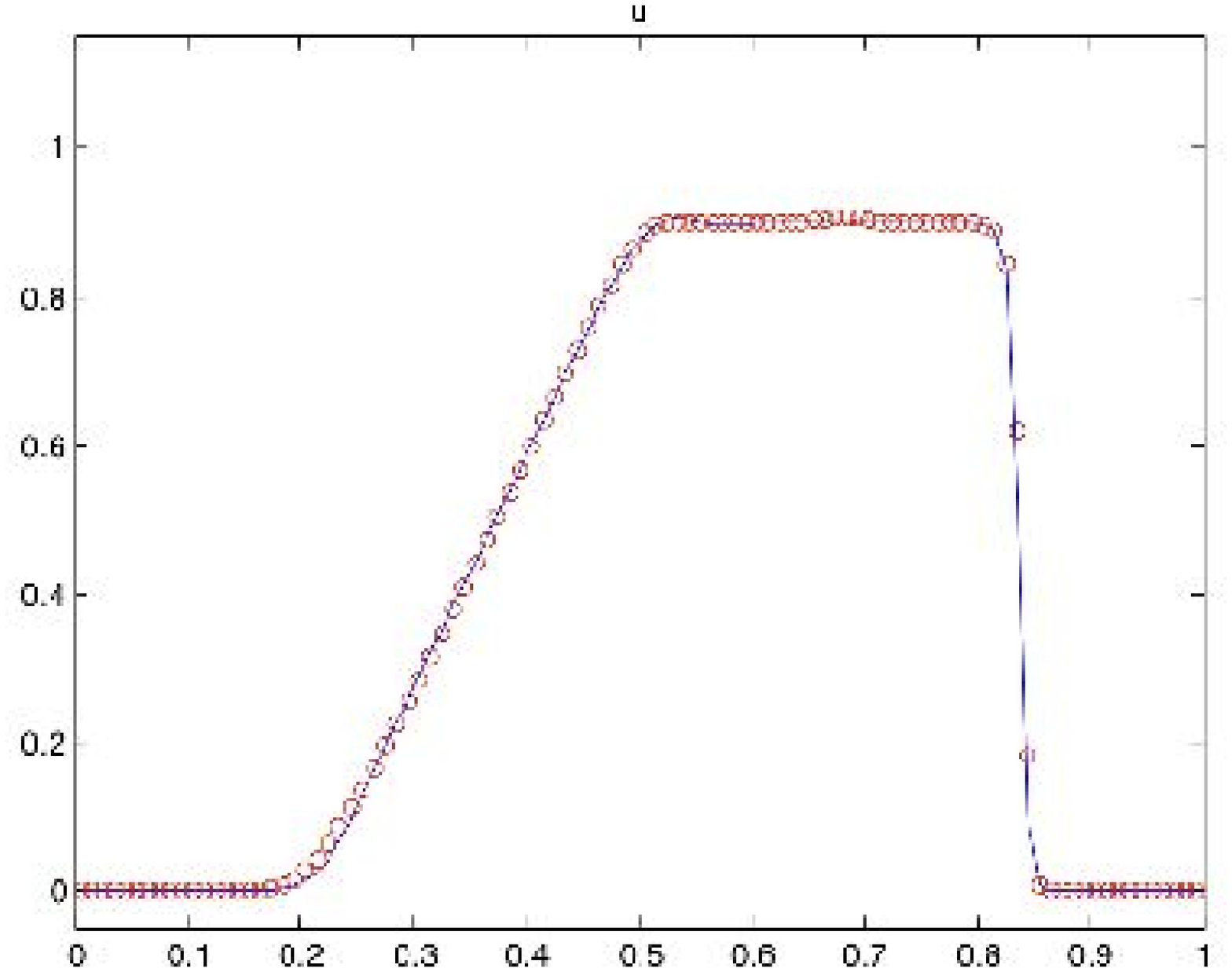}
\includegraphics[width=7.25cm]{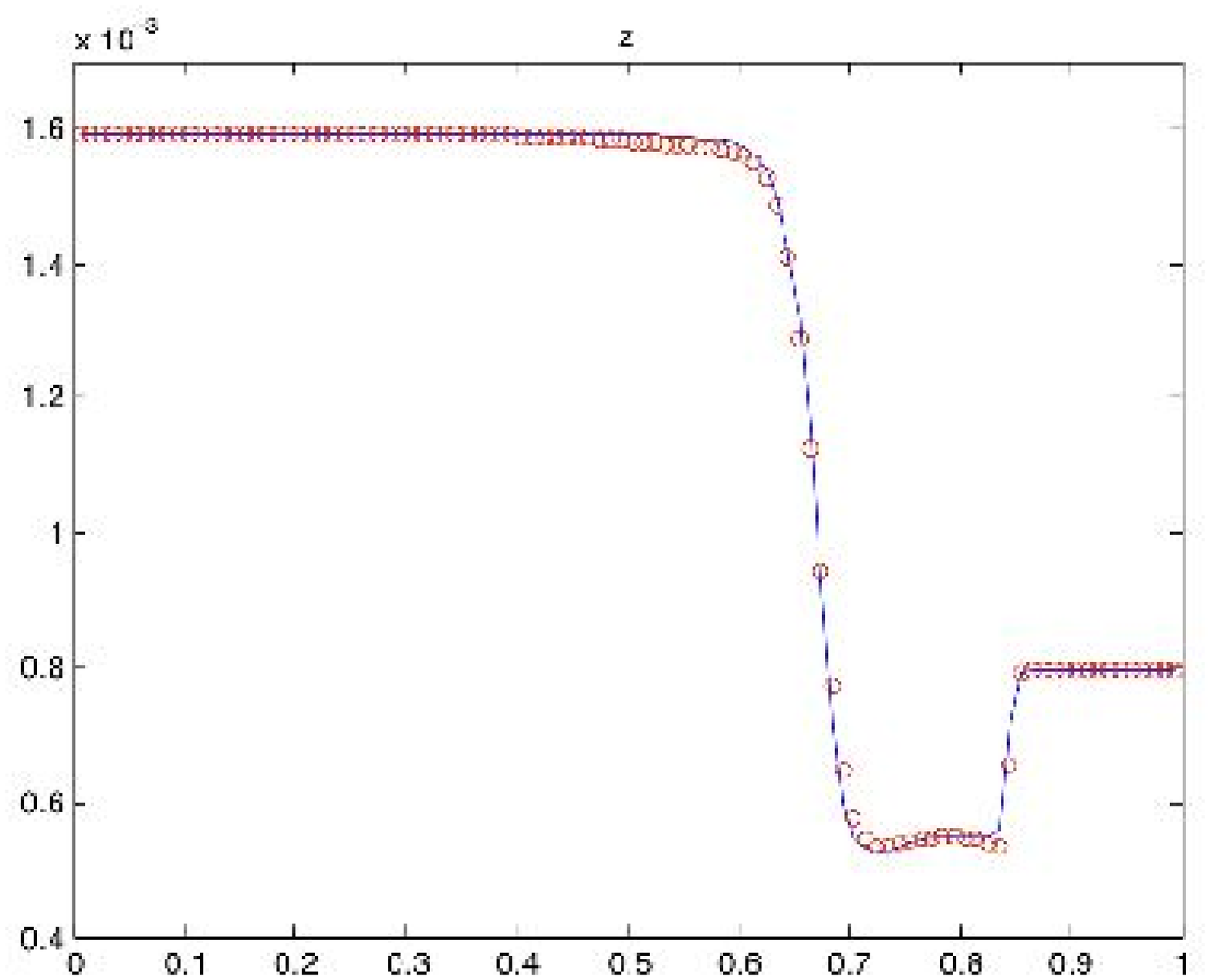}
\includegraphics[width=7.25cm]{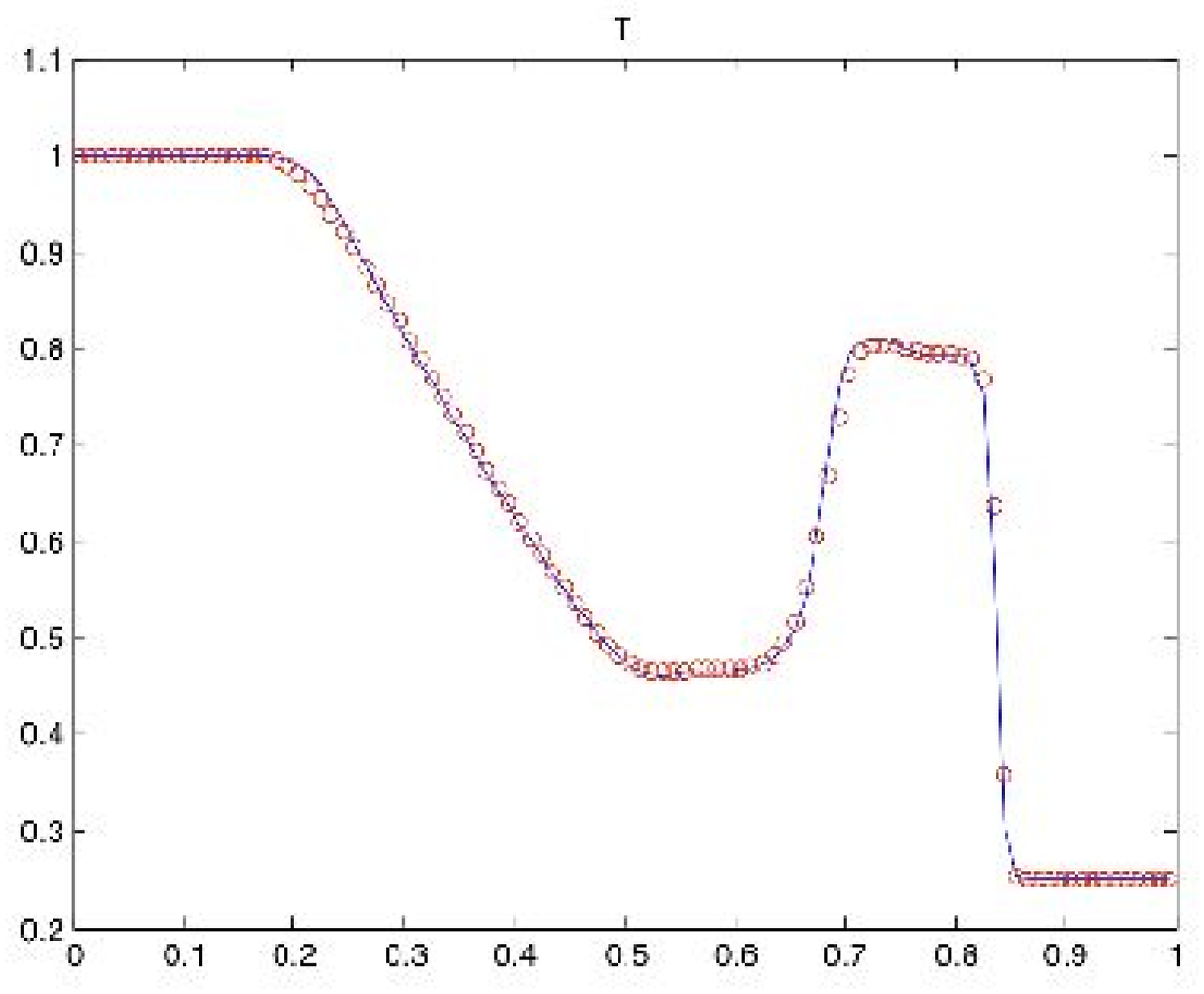}
\caption{Bose gas. $\epsilon=1e-4$, $\theta_0=0.01$, $z_l=0.0016$, $z_r=7.9546e-04$. Density $\rho$, velocity $u$, fugacity $z$ and temperature $T$ at $t=0.2$. $\Delta t=0.0013$, $\Delta x=0.01$. Solid line: KFVS scheme \cite{hu_KFVS} for quantum Euler equations (\ref{QEE}); $\circ$: New scheme (\ref{scheme1}) for QBE (\ref{QBE}).}
\label{Bose1e4_1}
\end{figure}

\begin{figure} [h!]
\centering
\includegraphics[width=7.25cm]{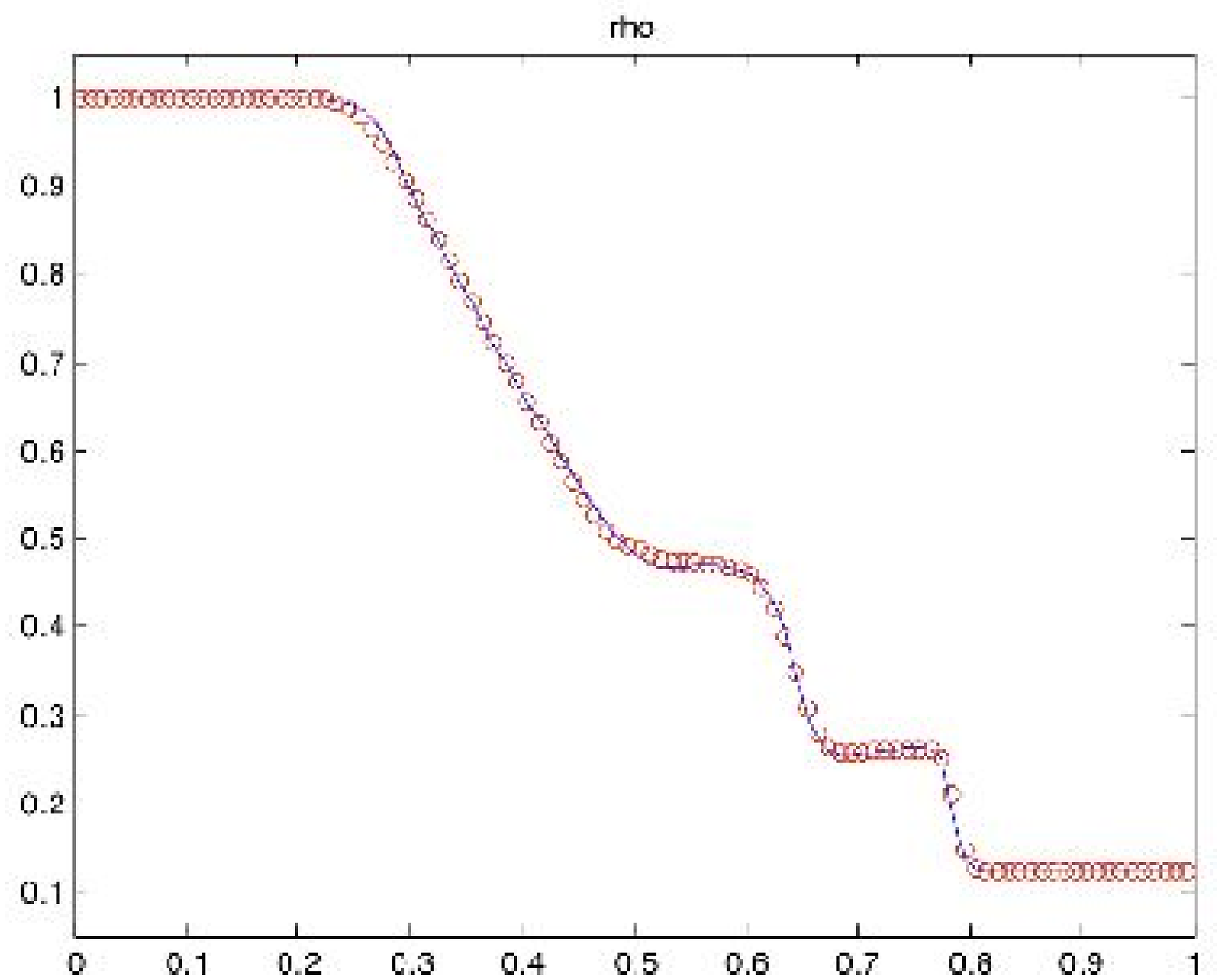}
\includegraphics[width=7.25cm]{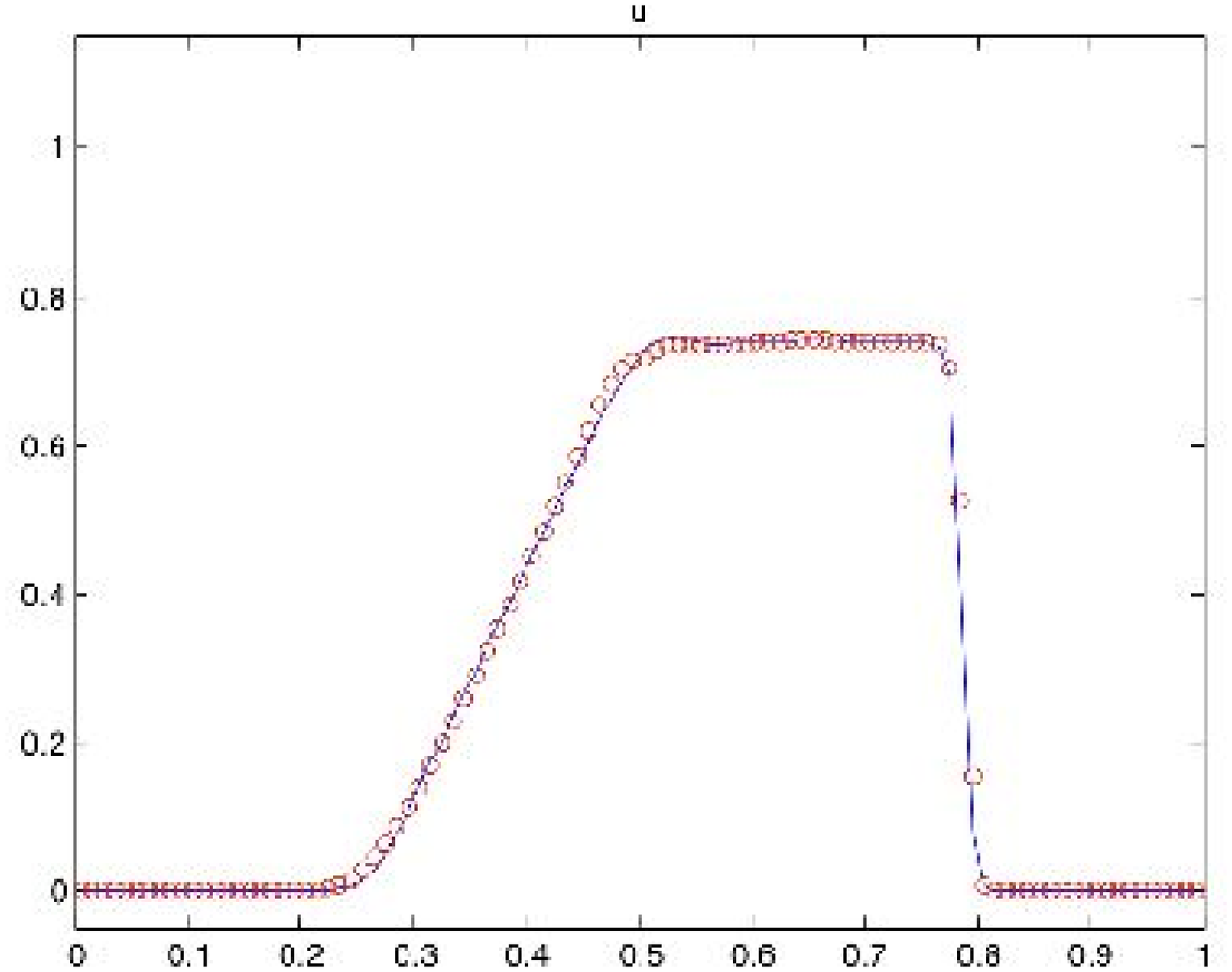}
\includegraphics[width=7.25cm]{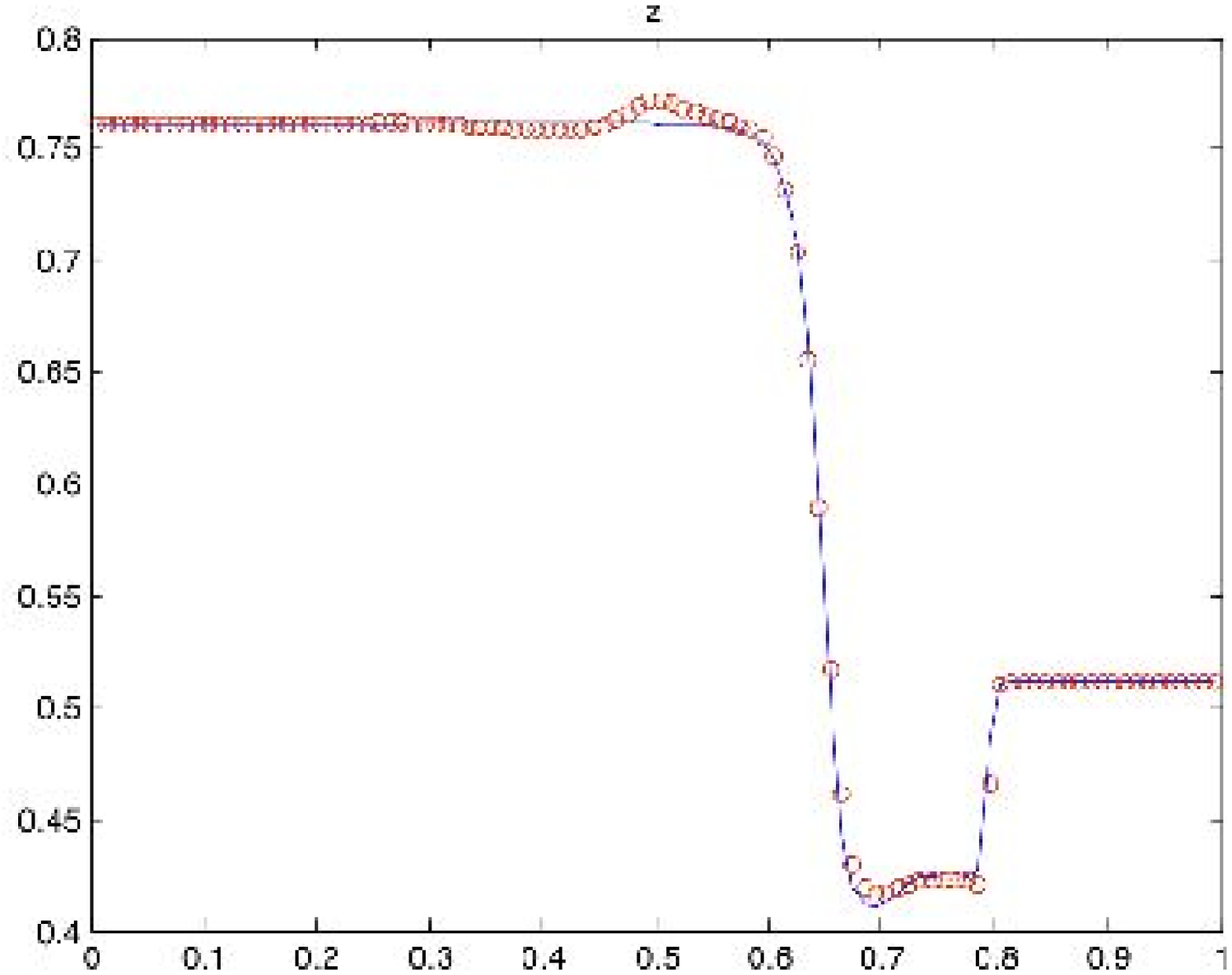}
\includegraphics[width=7.25cm]{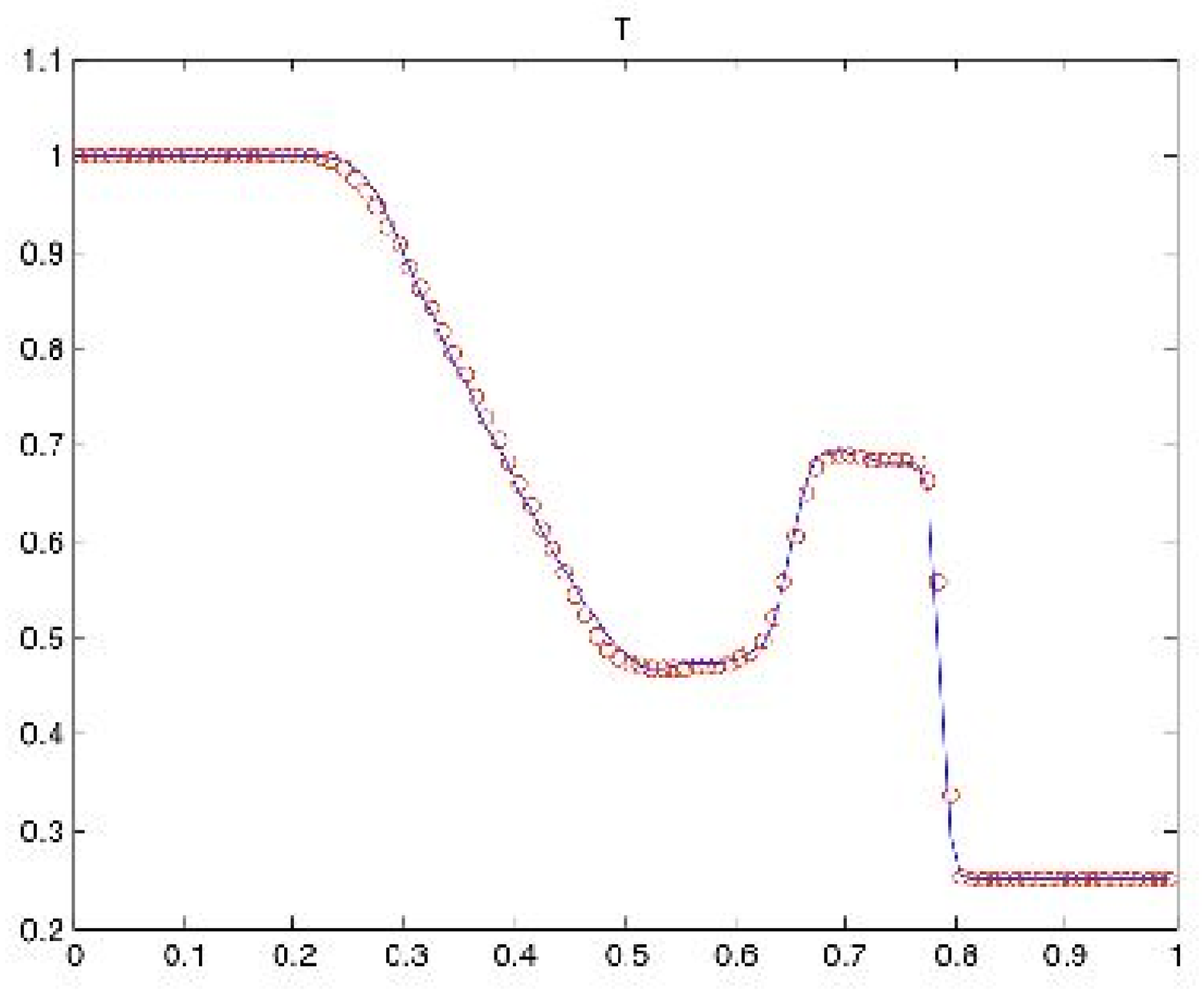}
\caption{Bose gas. $\epsilon=1e-4$, $\theta_0=9$, $z_l=0.7613$, $z_r=0.5114$. Density $\rho$, velocity $u$, fugacity $z$ and temperature $T$ at $t=0.2$. $\Delta t=0.0017$, $\Delta x=0.01$. Solid line: KFVS scheme \cite{hu_KFVS} for quantum Euler equations (\ref{QEE}); $\circ$: New scheme (\ref{scheme1}) for QBE (\ref{QBE}).}
\label{Bose1e4_3}
\end{figure}

\begin{figure} [h!]
\centering
\includegraphics[width=7.25cm]{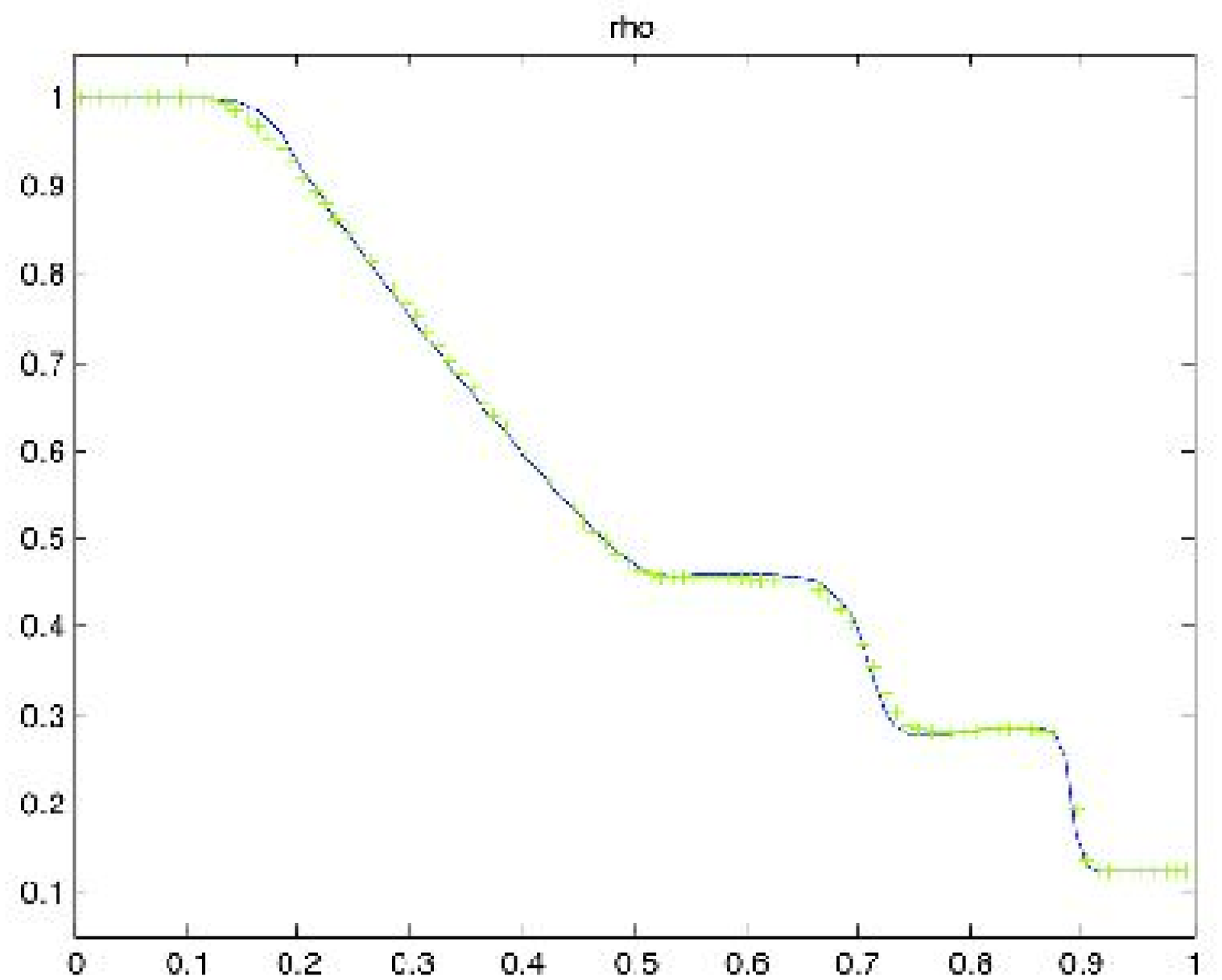}
\includegraphics[width=7.25cm]{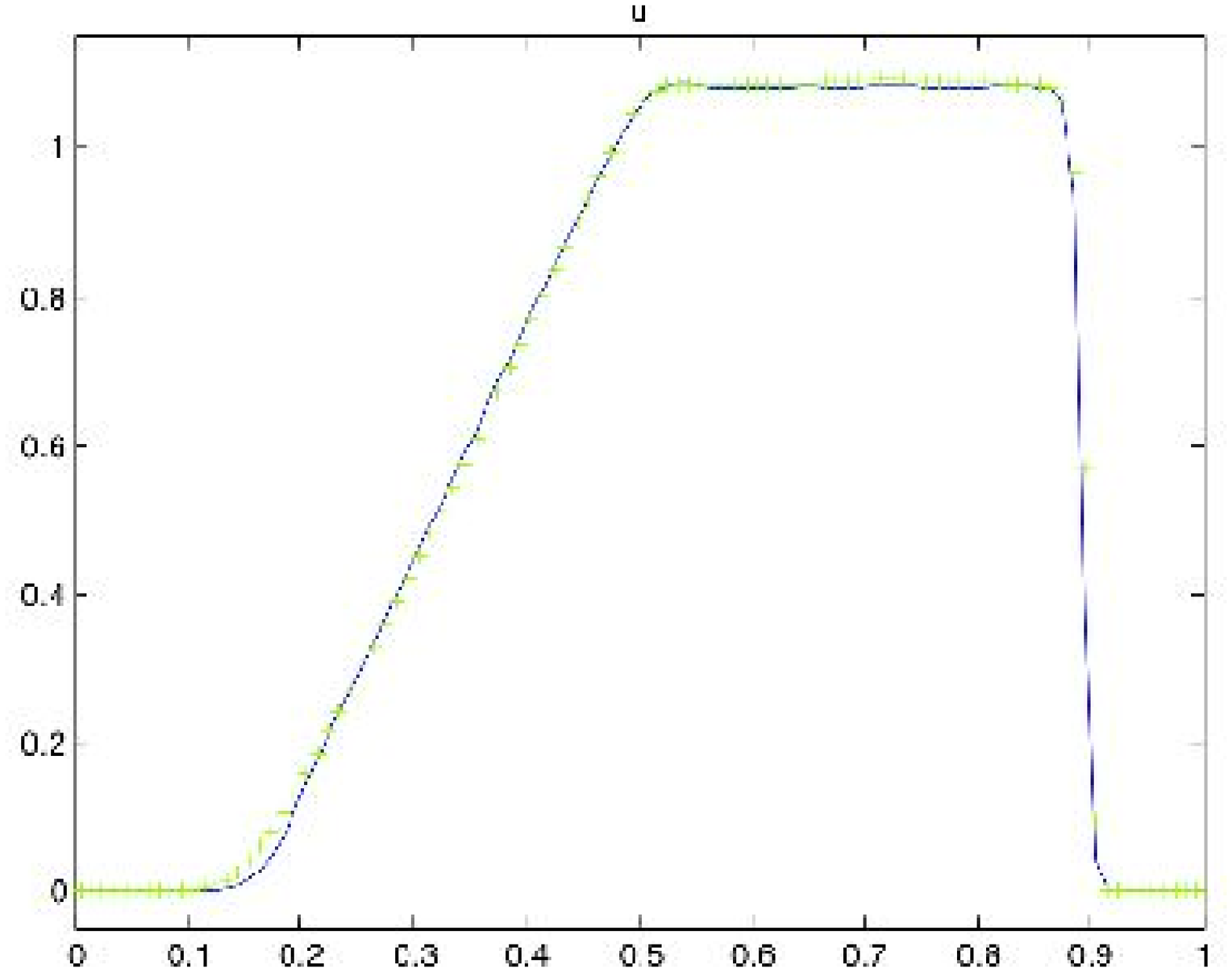}
\includegraphics[width=7.25cm]{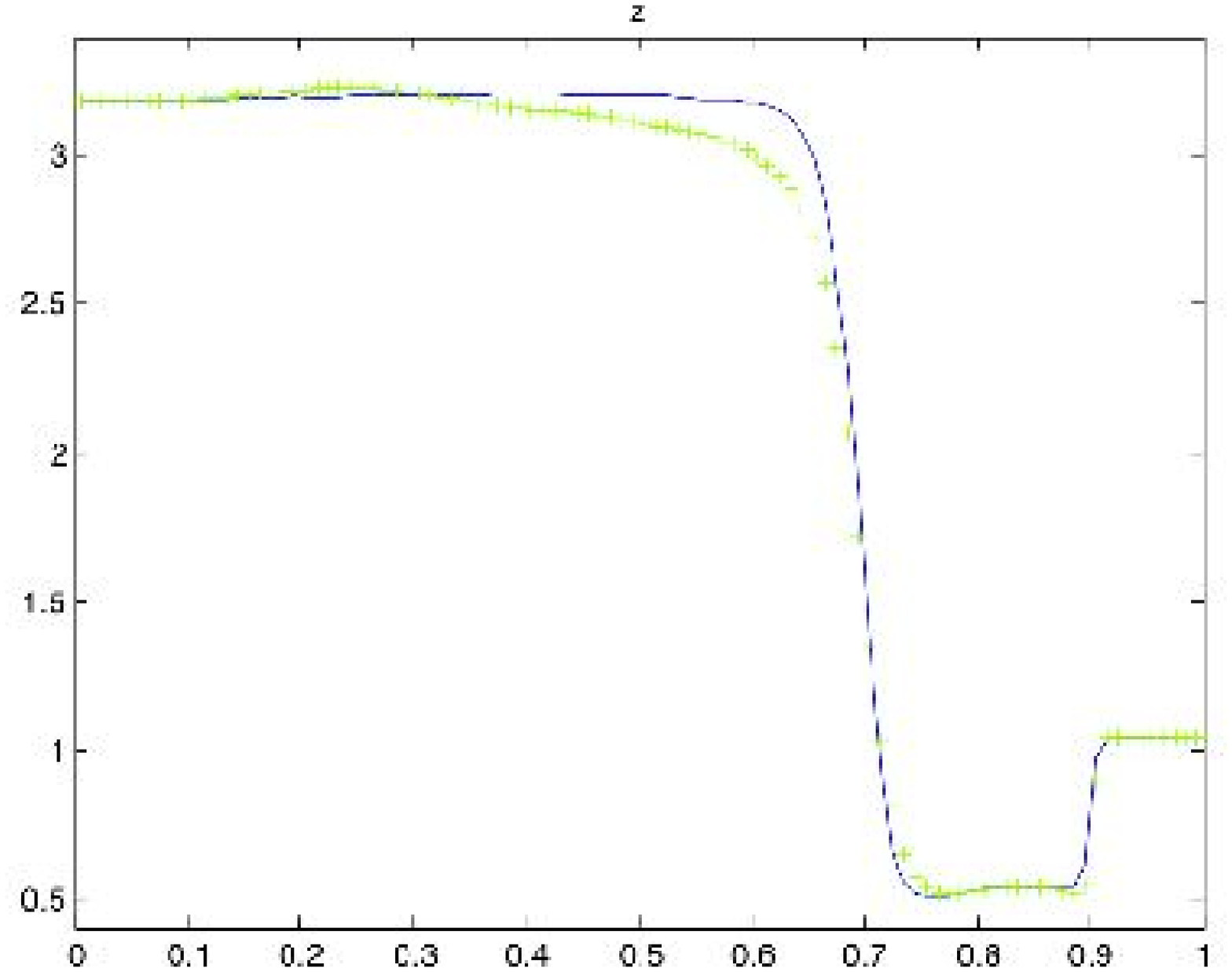}
\includegraphics[width=7.25cm]{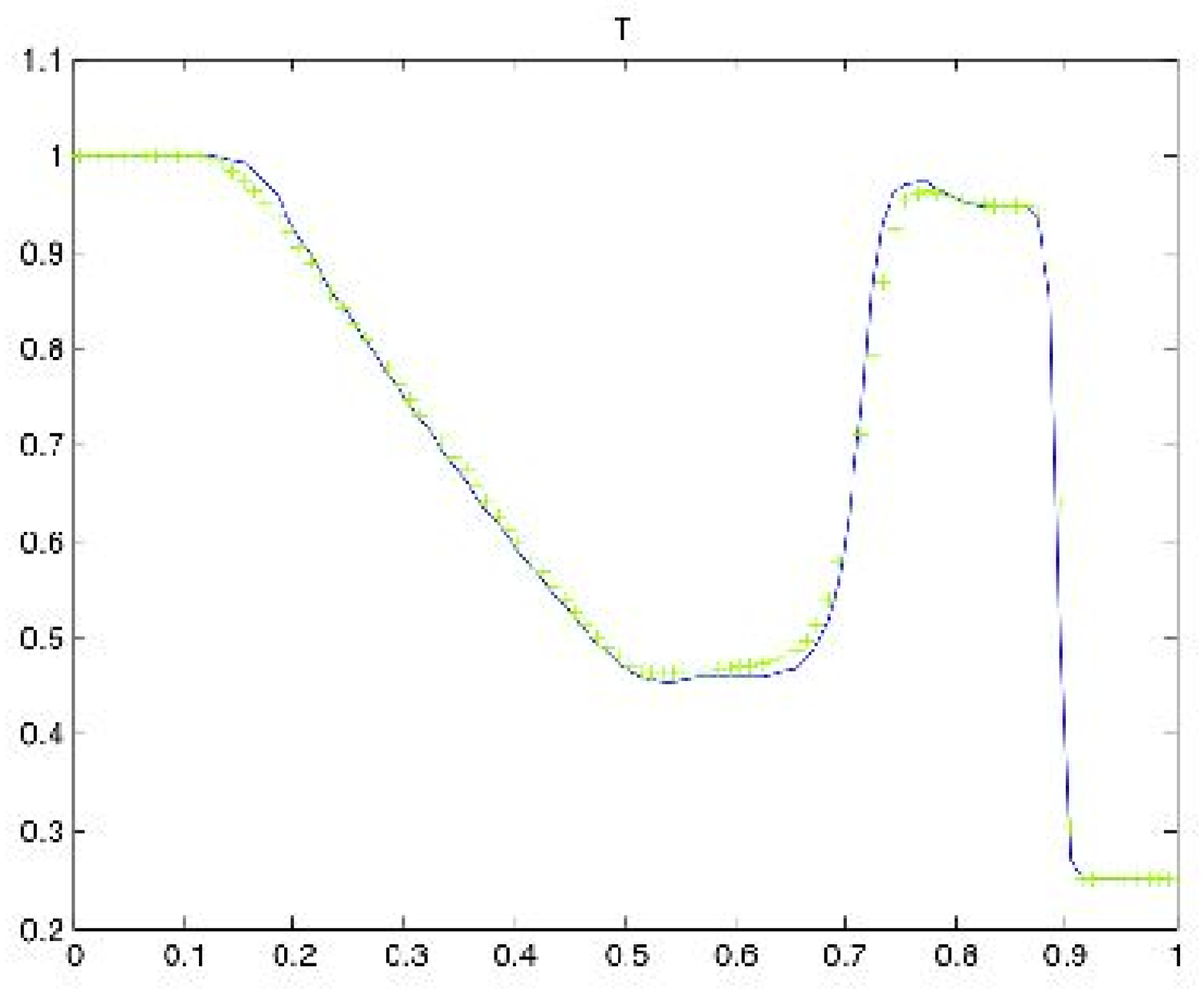}
\caption{Fermi gas. $\epsilon=1e-4$, $\theta_0=9$, $z_l=3.1887$, $z_r=1.0466$. Density $\rho$, velocity $u$, fugacity $z$ and temperature $T$ at $t=0.2$. $\Delta t=0.0013$, $\Delta x=0.01$. Solid line: KFVS scheme \cite{hu_KFVS} for quantum Euler equations (\ref{QEE}); $\circ$: New scheme (\ref{scheme1}) for QBE (\ref{QBE}).}
\label{Fermi1e4_3}
\end{figure}

We compare the results of our new scheme (\ref{scheme1}) with the kinetic scheme (KFVS scheme in \cite{hu_KFVS}) for the quantum Euler equations (\ref{QEE}). The time step $\Delta t$ is chosen by the CFL condition, independent of $\epsilon$. Fig.\ref{Bose1e4_1} shows the behaviors of a Bose gas when $\theta_0=0.01$. Fig.\ref{Bose1e4_3} shows the behaviors of a Bose gas when $\theta_0=9$. The solutions of a Fermi gas at $\theta_0=0.01$ are very similar to Fig.\ref{Bose1e4_1}, so we omit them here. Fig.\ref{Fermi1e4_3} shows the behaviors of a Fermi gas when $\theta_0=9$. All the results agree well in this regime, which exactly implies the scheme (\ref{scheme1}) is asymptotic preserving (when the Knudsen number $\epsilon$ goes to zero, the scheme becomes a fluid solver).

\subsection{Kinetic Regime}
\begin{figure} [h!]
\centering
\includegraphics[width=7.25cm]{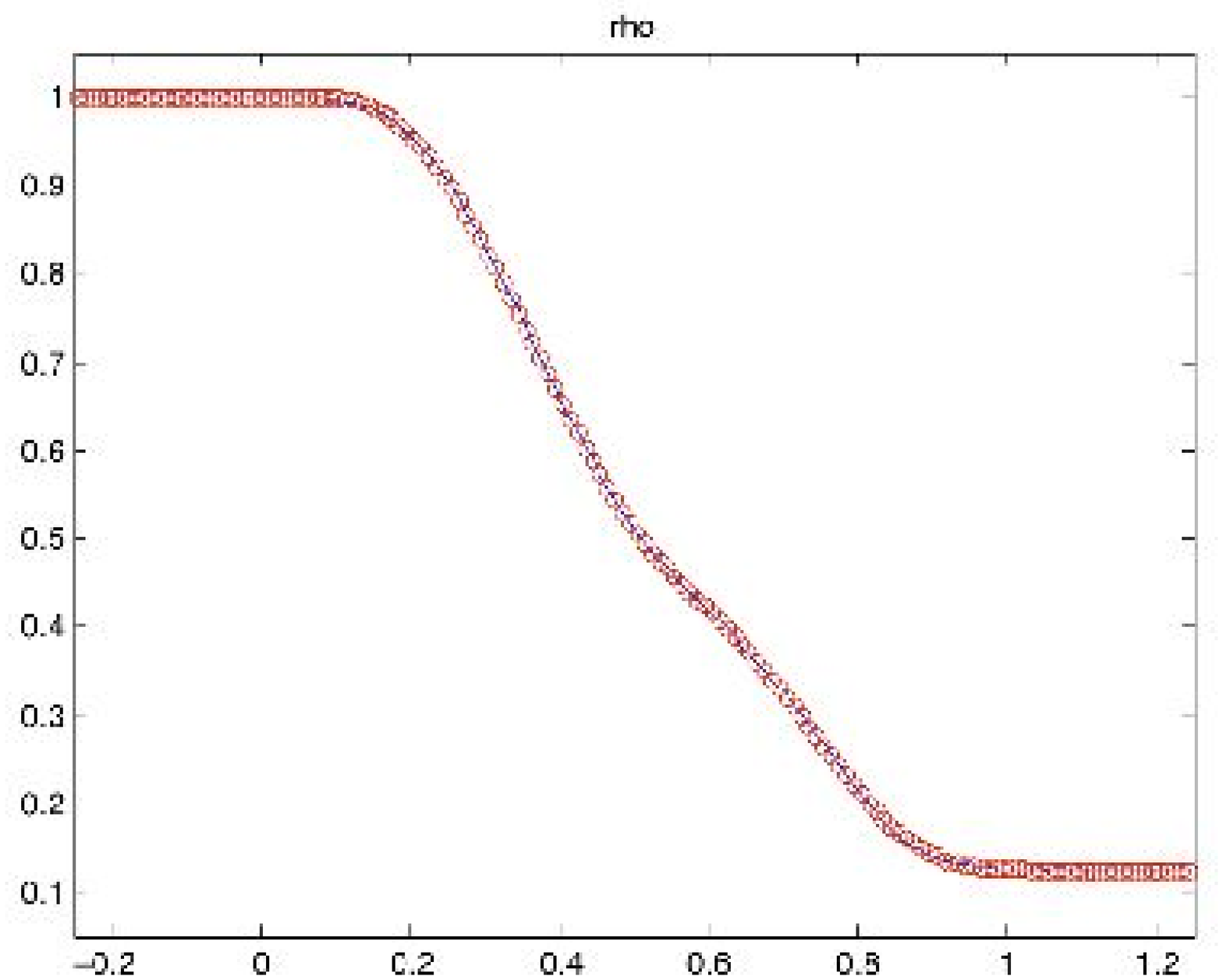}
\includegraphics[width=7.25cm]{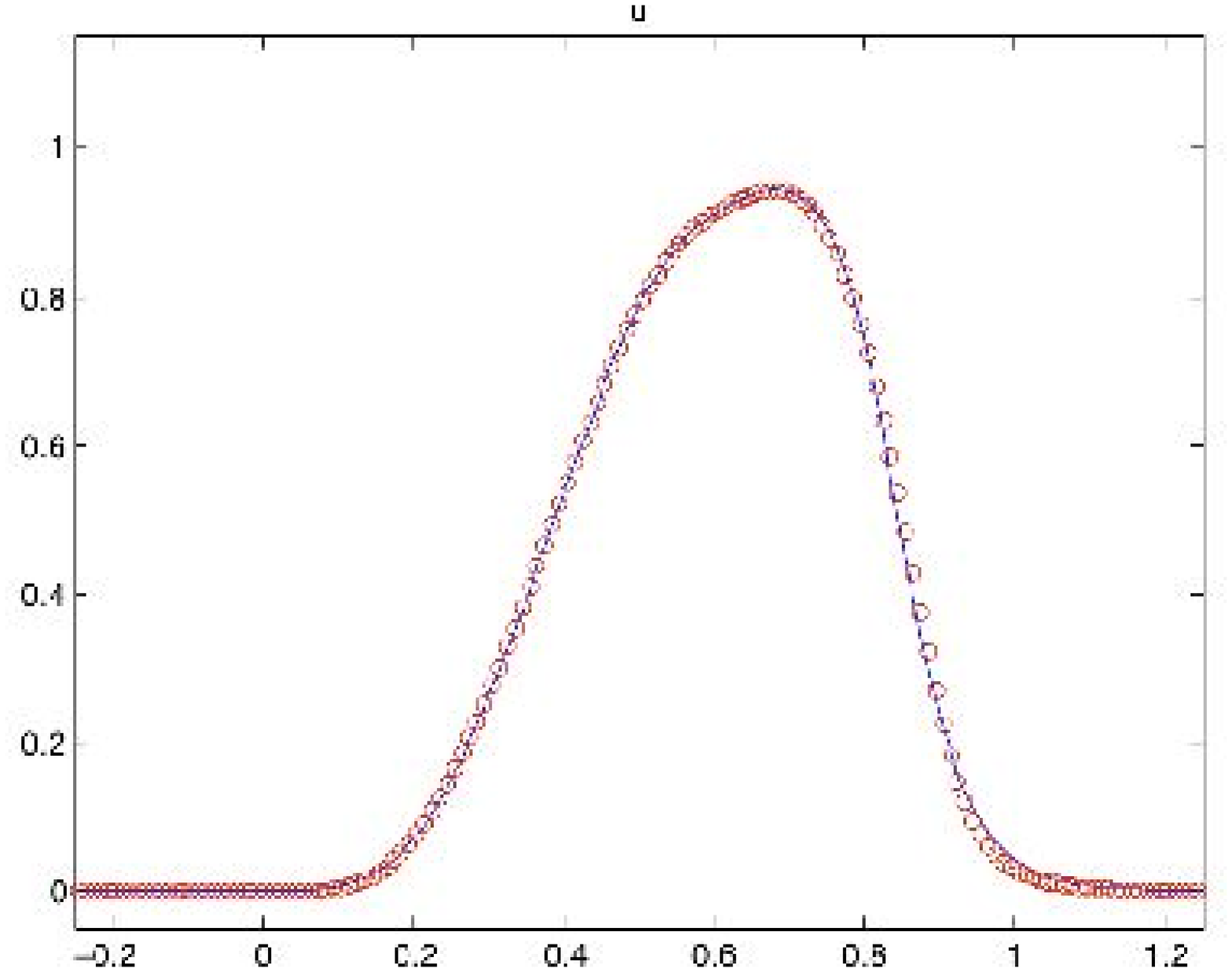}
\includegraphics[width=7.25cm]{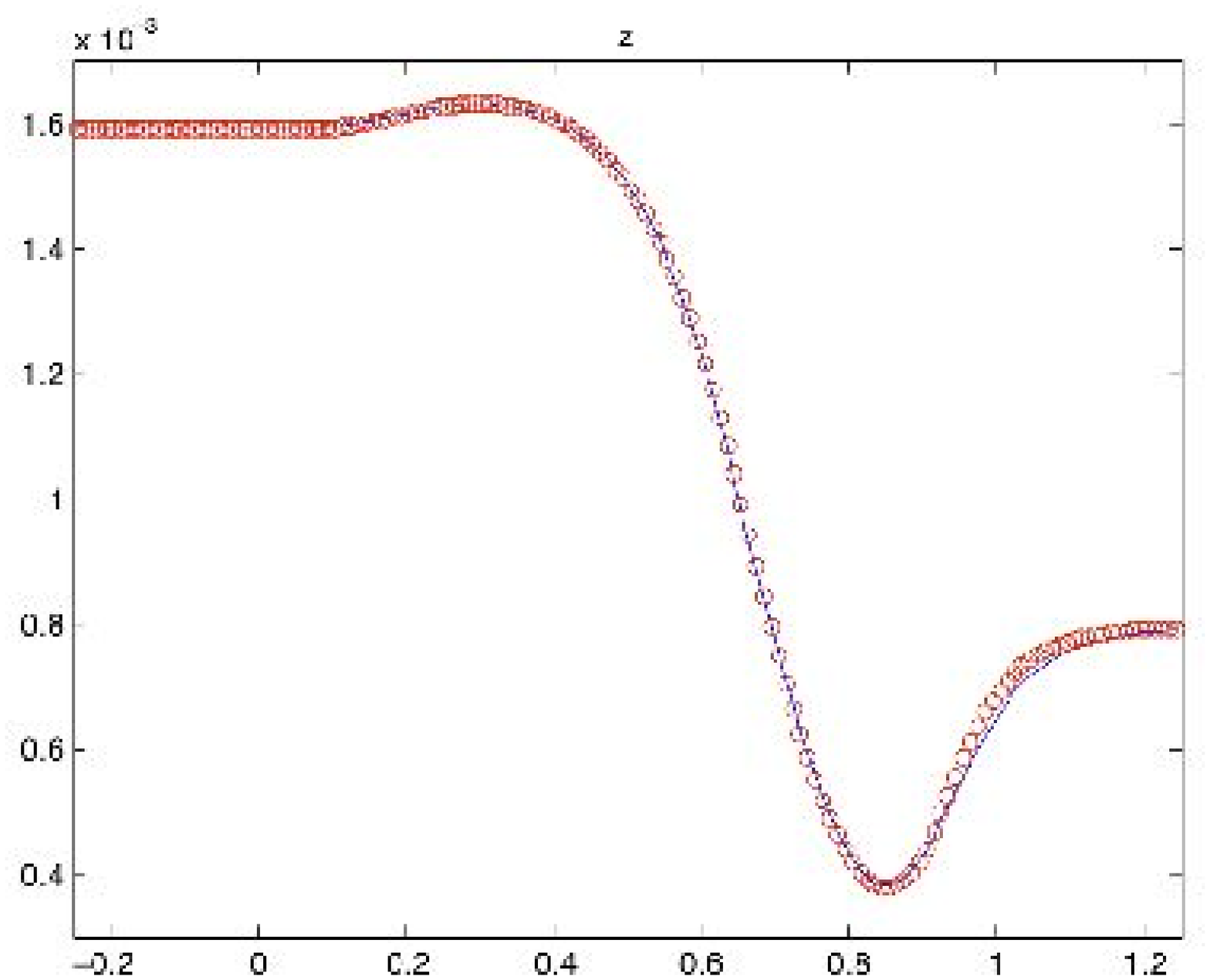}
\includegraphics[width=7.25cm]{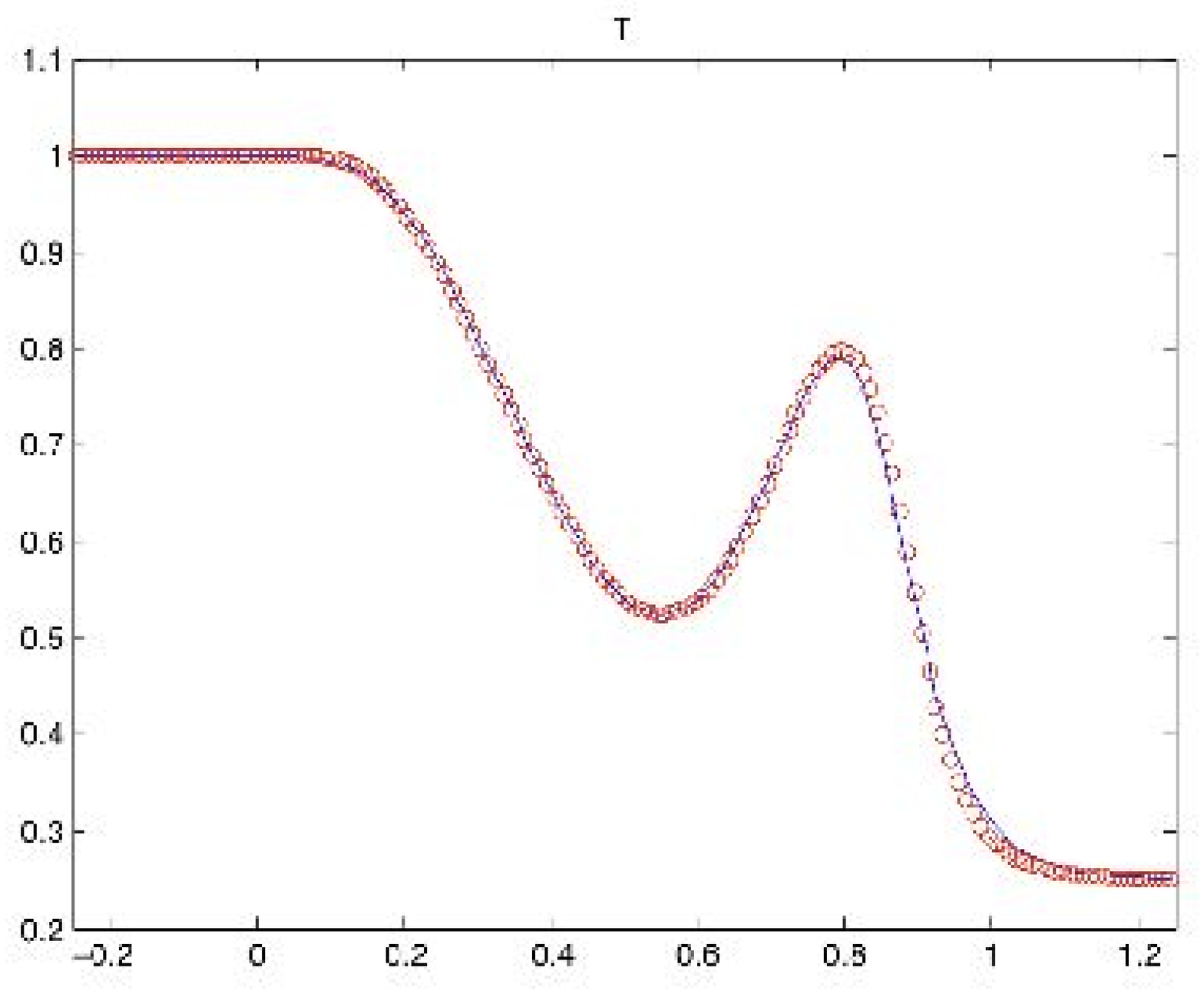}
\caption{Bose gas. $\epsilon=1e-2$, $\theta_0=0.01$, $z_l=0.0016$, $z_r=7.9546e-04$. Density $\rho$, velocity $u$, fugacity $z$ and temperature $T$ at $t=0.2$. $\Delta t=0.0013$, $\Delta x=0.01$. Solid line: Forward Euler scheme for QBE (\ref{QBE}); $\circ$: New scheme (\ref{scheme1}) for QBE (\ref{QBE}).}
\label{Bose1e2_1}
\end{figure}

\begin{figure} [h!]
\centering
\includegraphics[width=7.25cm]{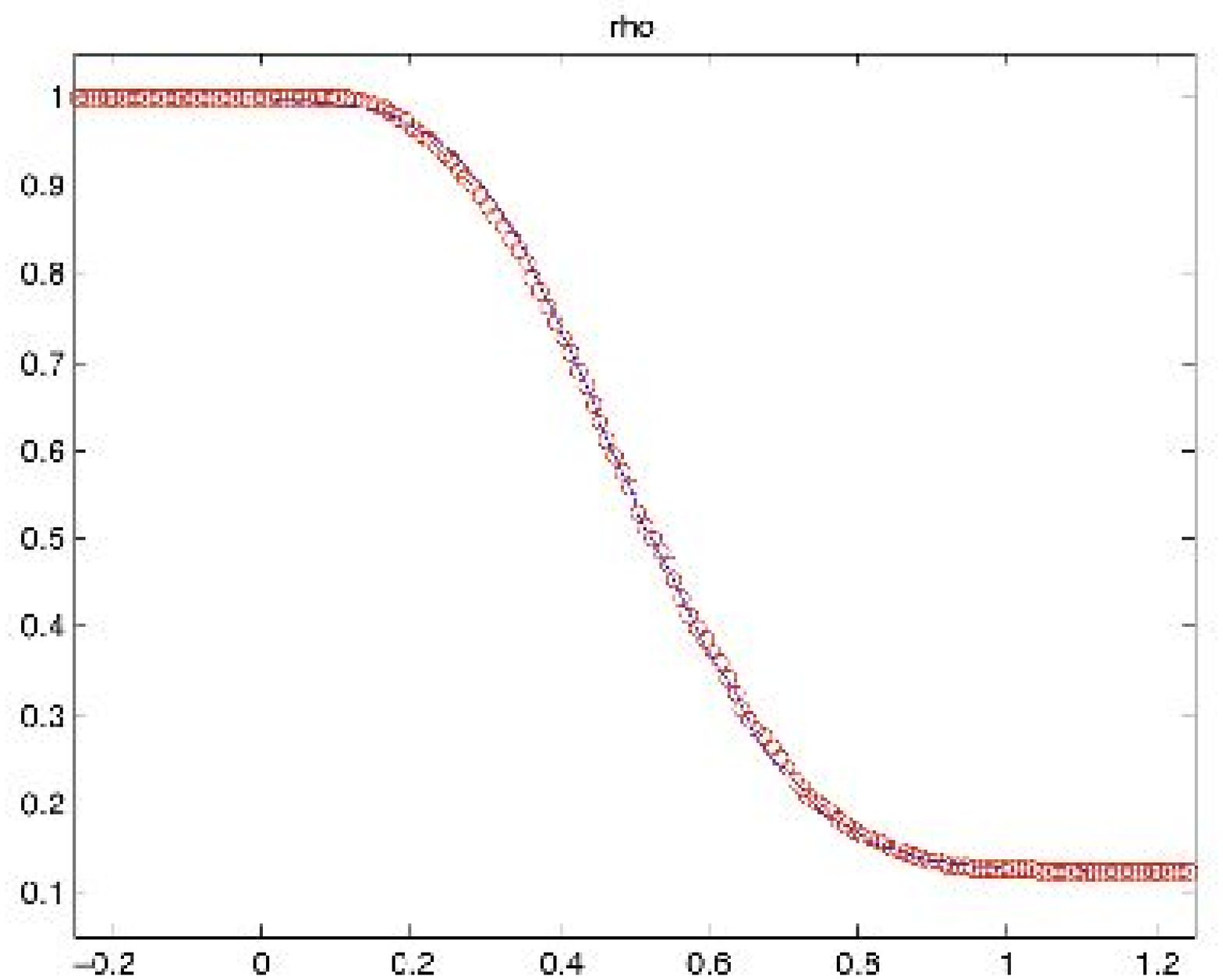}
\includegraphics[width=7.25cm]{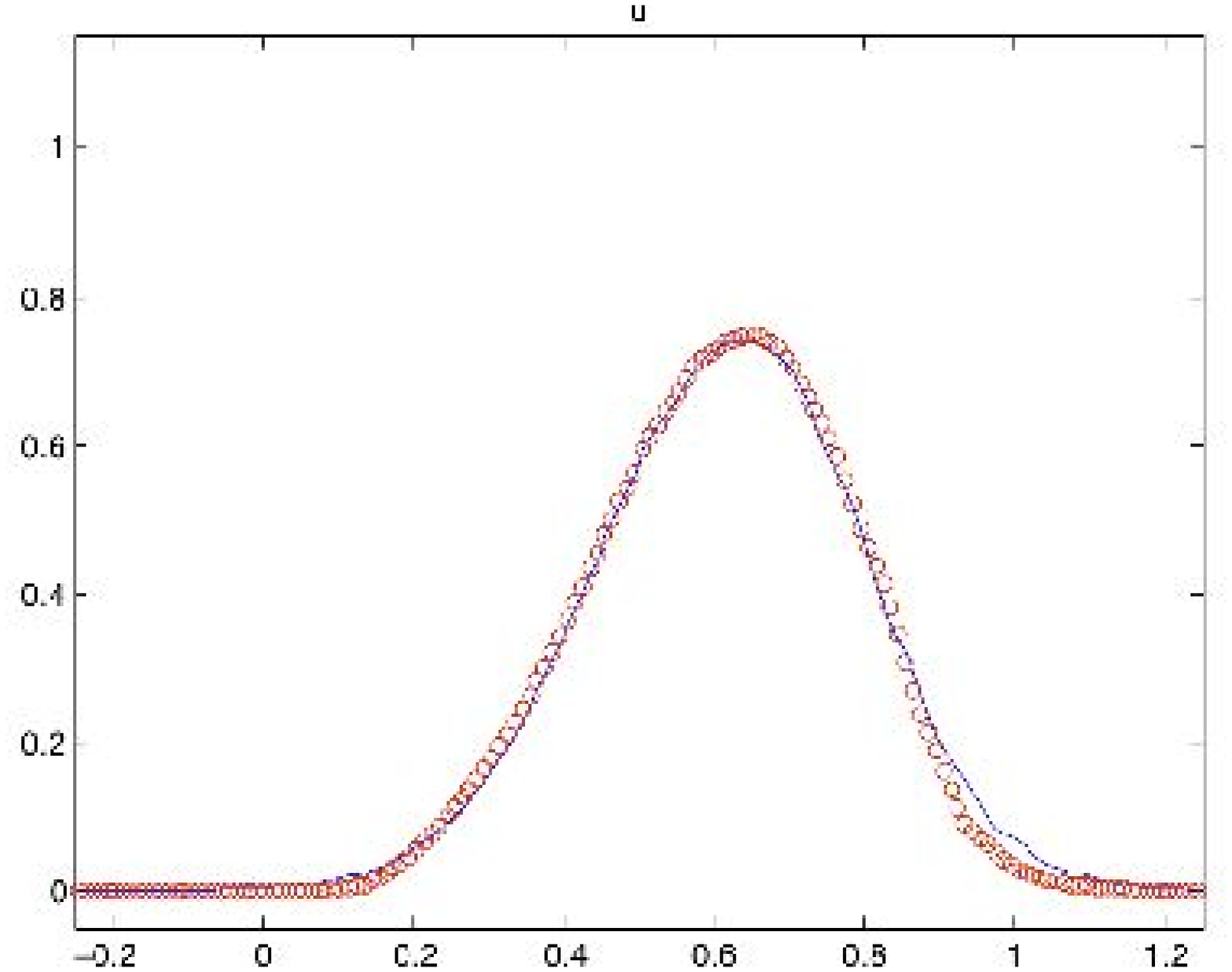}
\includegraphics[width=7.25cm]{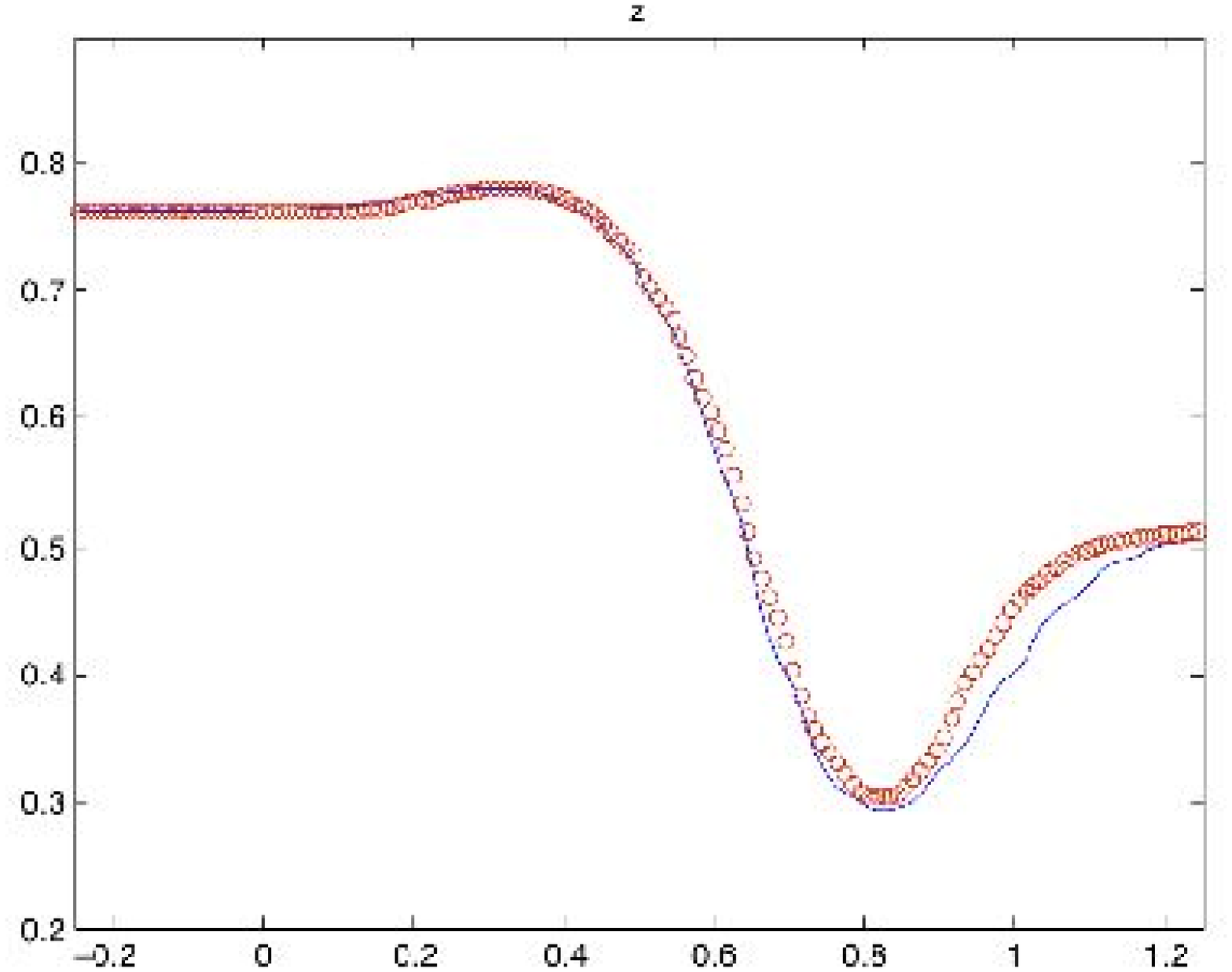}
\includegraphics[width=7.25cm]{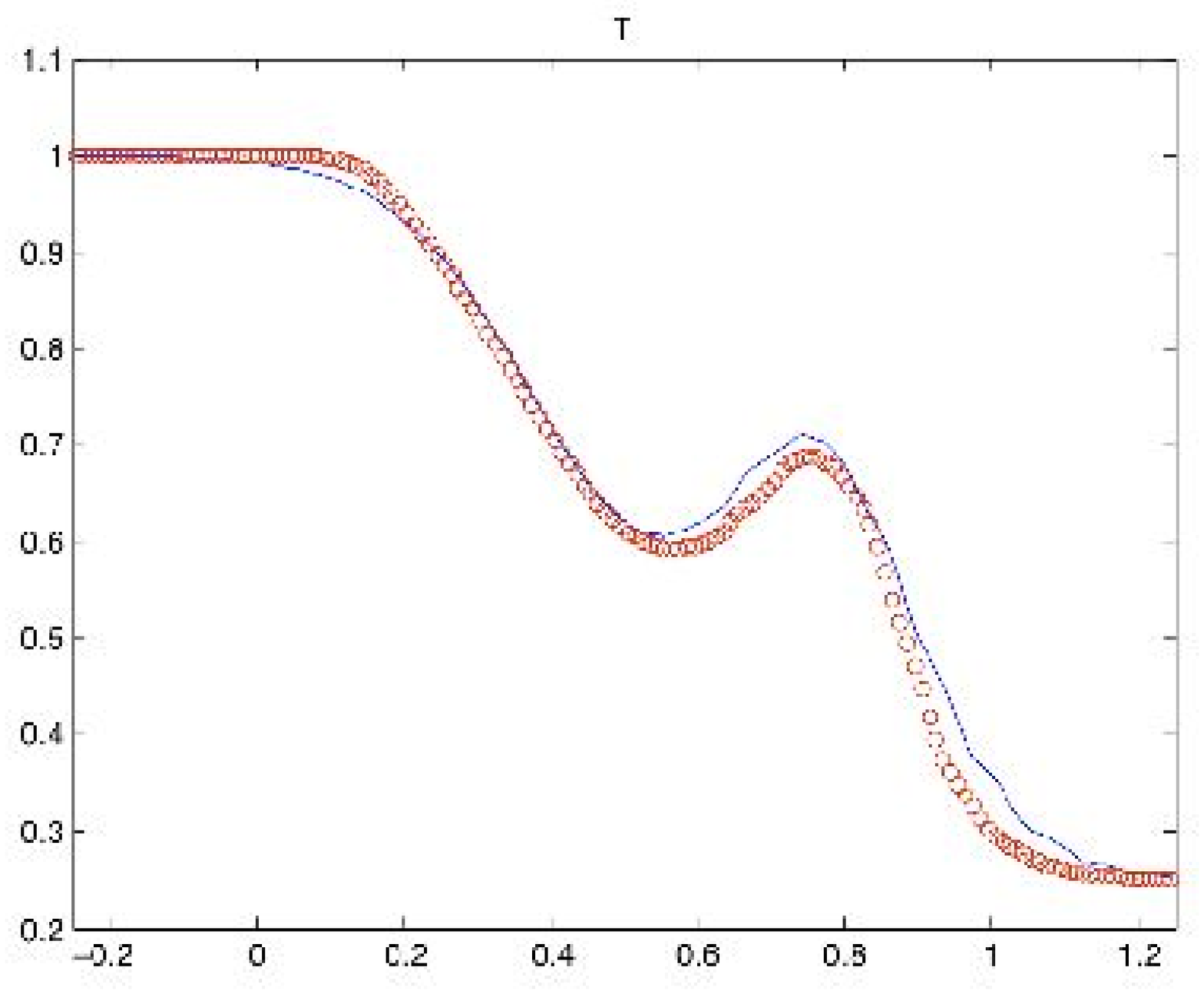}
\caption{Bose gas. $\epsilon=1e-1$, $\theta_0=9$, $z_l=0.7613$, $z_r=0.5114$. Density $\rho$, velocity $u$, fugacity $z$ and temperature $T$ at $t=0.2$. $\Delta t=0.0017$, $\Delta x=0.01$. Solid line: Forward Euler scheme for QBE (\ref{QBE}); $\circ$: New scheme (\ref{scheme1}) for QBE (\ref{QBE}).}
\label{Bose1e1_3}
\end{figure}

\begin{figure} [h!]
\centering
\includegraphics[width=7.25cm]{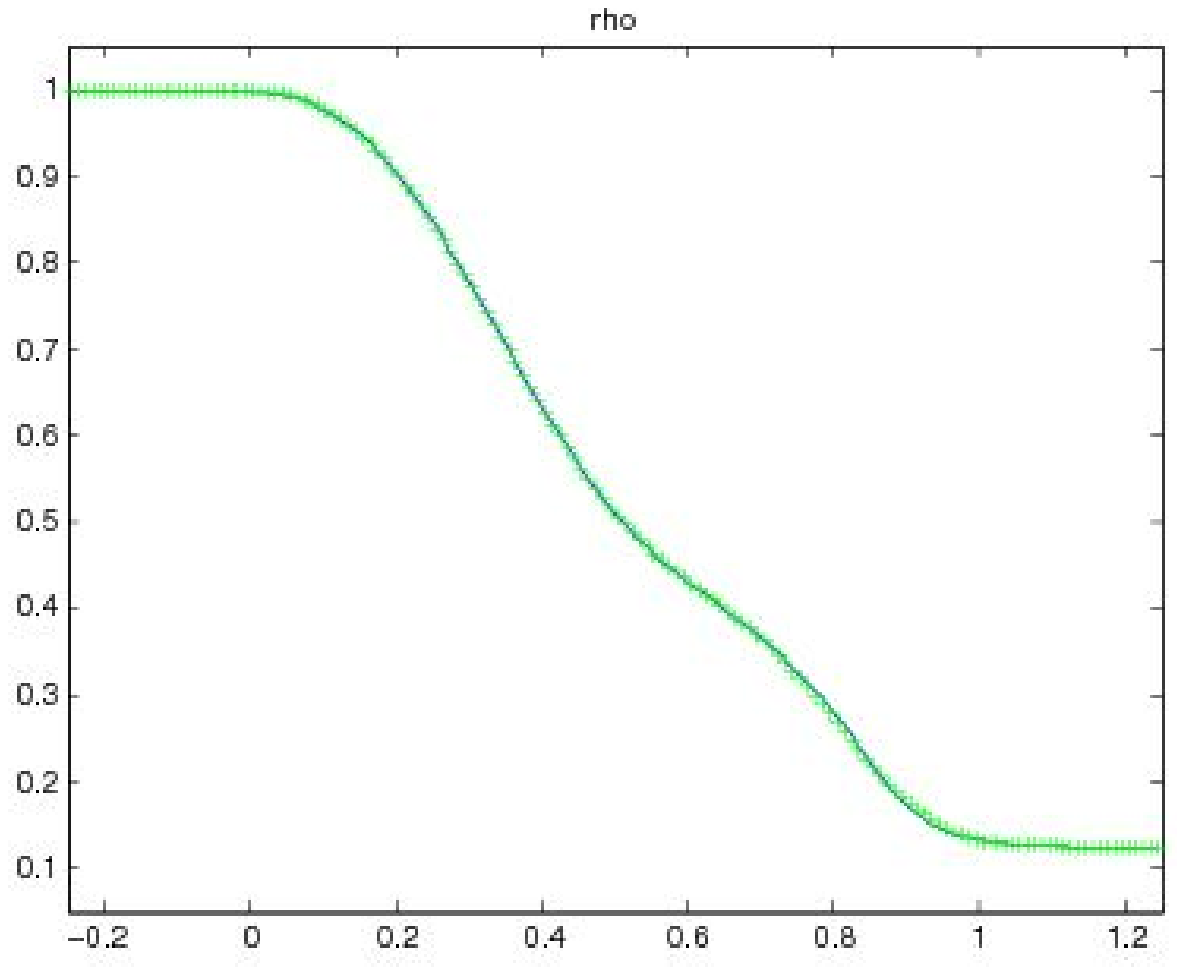}
\includegraphics[width=7.25cm]{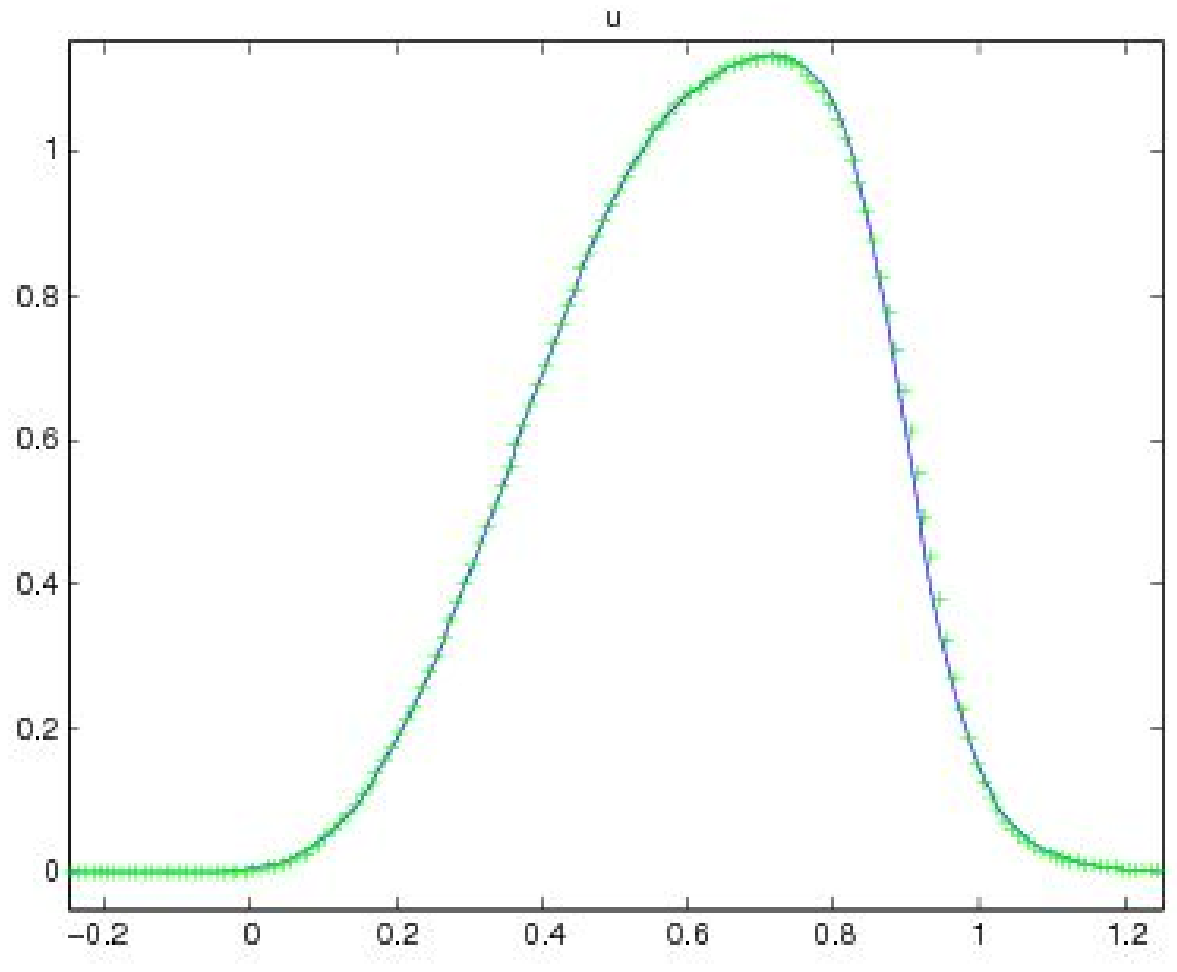}
\includegraphics[width=7.25cm]{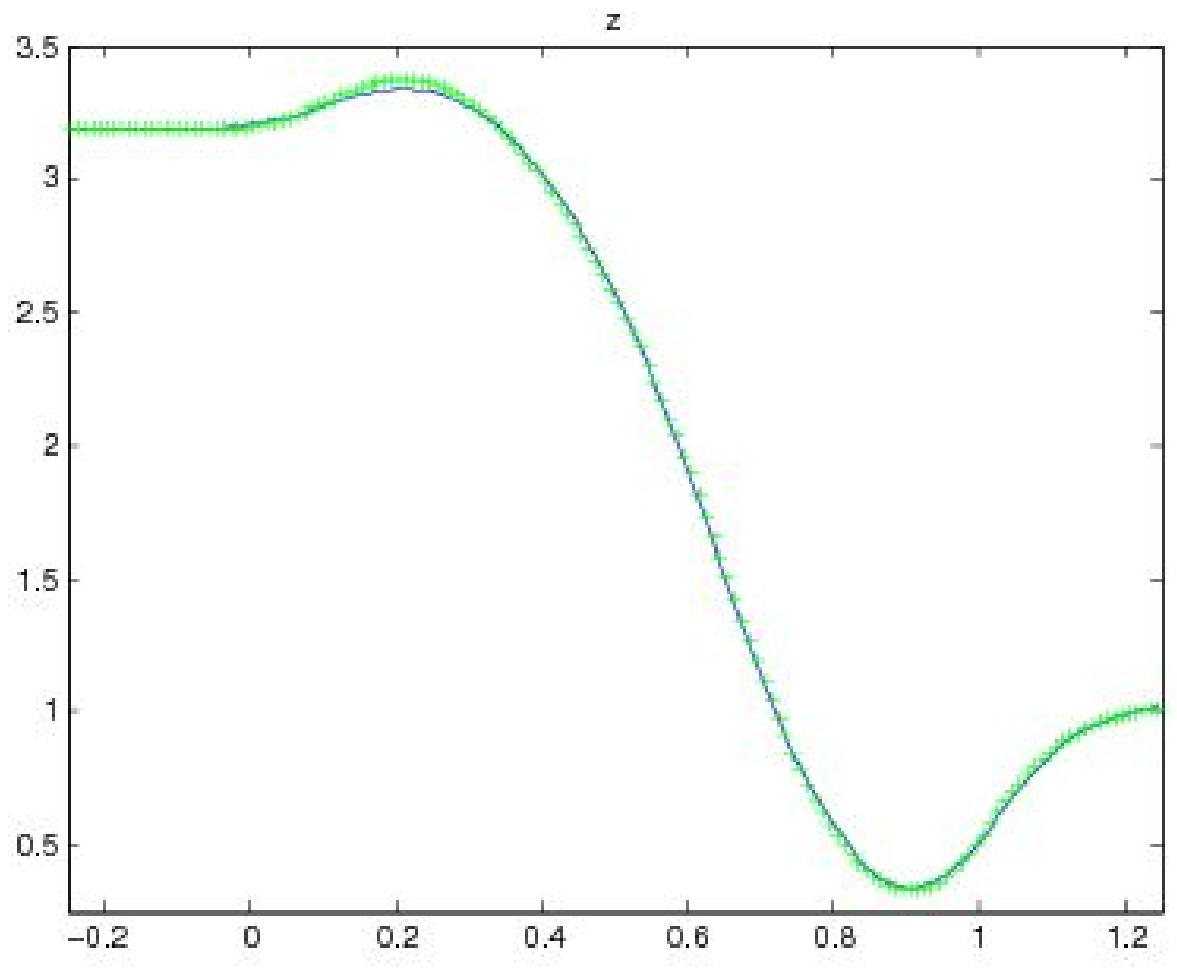}
\includegraphics[width=7.25cm]{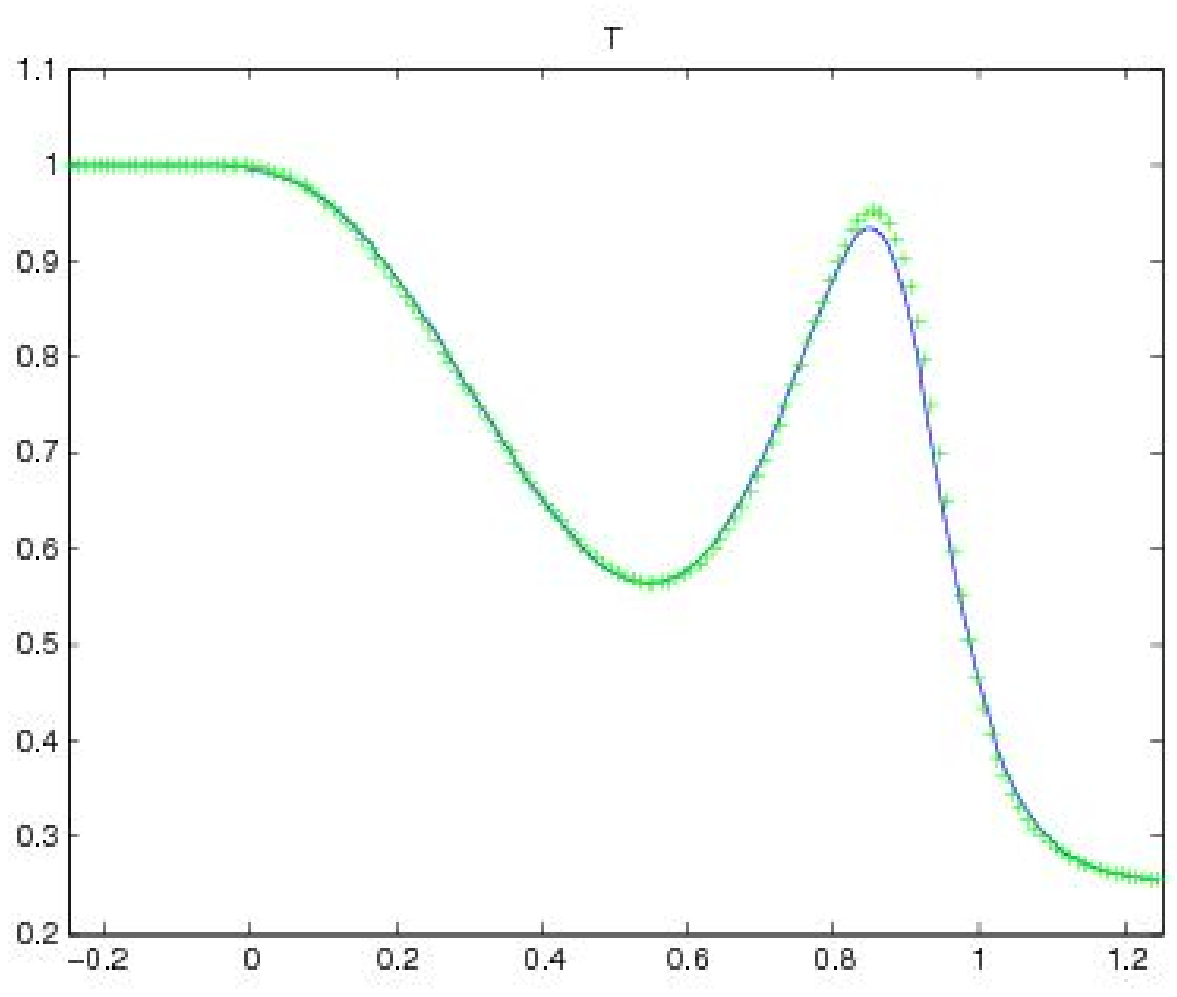}
\caption{Fermi gas. $\epsilon=1e-2$, $\theta_0=9$, $z_l=3.1887$, $z_r=1.0466$. Density $\rho$, velocity $u$, fugacity $z$ and temperature $T$ at $t=0.2$. $\Delta t=0.0013$, $\Delta x=0.01$. Solid line: Forward Euler scheme for QBE (\ref{QBE}); $\circ$: New scheme (\ref{scheme1}) for QBE (\ref{QBE}).}
\label{Fermi1e2_3}
\end{figure}

We compare the results of our new scheme (\ref{scheme1}) with the explicit forward Euler scheme. The time step $\Delta t$ for the new scheme is still chosen by the CFL condition. When the Knudsen number $\epsilon$ is not very small, $10^{-1}$ or $10^{-2}$, the above $\Delta t$ is also enough for the explicit scheme. Fig.\ref{Bose1e2_1} shows the behaviors of a Bose gas when $\theta_0=0.01$. Fig.\ref{Bose1e1_3} shows the behaviors of a Bose gas when $\theta_0=9$. The solutions of a Fermi gas at $\theta_0=0.01$ are very similar to Fig.\ref{Bose1e2_1}, so we omit them here. Fig.\ref{Fermi1e2_3} shows the behaviors of a Fermi gas when $\theta_0=9$. Again all the results agree well which means the scheme (\ref{scheme1}) is also reliable in the kinetic regime. To avoid the boundary effect, all the simulations in this subsection were carried out on a slightly larger spatial domain $x\in[-0.25, 1.25]$.

\section{Conclusion}
A novel scheme was introduced for the quantum Boltzmann equation, starting from the scheme in \cite{FJ10}. The new idea here is to penalize the quantum collision operator by a `classical' BGK operator so as to avoid the difficulty of inverting the nonlinear system $\rho=\rho(z,T)$, $e=e(z,T)$. The new scheme is uniformly stable in terms of the Knudsen number, and can capture the fluid (Euler) limit even if the small scale is not numerically resolved. We have also developed a spectral method for the quantum collision operator, following its classical counterpart \cite{MP06,FMP06}.

So far we have not considered the quantum gas in the extreme case. For example, the Bose gas becomes degenerate when the fugacity $z=1$. Many interesting phenomena happen in this regime. Our future work will focus on this aspect.

\vspace{5mm}

\noindent \textbf{Acknowledgments.} The second author would like to thank Mr. Bokai Yan for helpful discussions on the spectral method of the collision operator.


\begin{flushleft} 
\signff 
\end{flushleft}
\vspace{-4.25cm}
\begin{flushright} 
\signjh 
\end{flushright}

\begin{flushleft} 
\signsj
\end{flushleft}

\end{document}